\documentclass[12pt]{article}
\usepackage{amssymb}
\usepackage{tikz}
\usepackage{mathrsfs}
\font\eufm=eufm10 at 14pt\font\eufms=eufm10\font\eufmss=eufm7\newfam\eufam
\textfont\eufam=\eufm\scriptfont\eufam=\eufms\scriptscriptfont\eufam=\eufmss

\def\build#1_#2^#3{\mathrel{\mathop{\kern 0pt#1}\limits_{#2}^{#3}}}
\def\Z{{\bf Z}}\def\Q{{\bf Q}}\def\C{{\bf C}}\def\F{{\bf F}}\def\H{{\rm H}}\def\Hom{{\rm Hom}}

\def\rde{\mathscr}

\def\cqfd{\hfill\vbox{\hrule\hbox{\vrule height6pt depth0pt\hskip 6pt \vrule height6pt}\hrule\relax}}
\def\noi{\noindent}
\def\hfl#1#2{\smash{\mathop{\hbox to 12 mm{\rightarrowfill}}\limits^{\scriptstyle#1}_{\scriptstyle#2}}}
\def\vfl#1#2{\llap{$\scriptstyle #1$}\left\downarrow\vbox to 6mm{}\right.\rlap{$\scriptstyle #2$}}
\def\diagram#1{\def\normalbaselines{\baselineskip=0pt\lineskip=10pt\lineskiplimit=1pt} \matrix{#1}}
\def\pv{\raise 2pt\hbox{$\bigwedge$}}
\begin{document}

\overfullrule=0pt
\
\vskip 64pt
\centerline{\bf Cohn localization of finite group rings.}
\vskip 12pt
\centerline{Pierre Vogel\footnote{Universit\'e Paris Diderot, Institut de Math\'ematiques de Jussieu-Paris Rive Gauche (UMR 7586),
B\^atiment Sophie Germain, Case 7012, 75205--Paris Cedex 13 --- Email: pierre.vogel@imj-prg.fr}}
\vskip 48pt
\noi{\bf Abstract.} The purpose of this paper is to give a complete description of the Cohn localization of the augmentation map $\Z[G]\rightarrow\Z$ when $G$ is any
finite group.
\vskip 12pt
\noi{\bf Keywords:} Non-commutative localization, category of chain complexes.
\vskip 24pt
\noi{\bf Introduction.} 
Let $W$ be a set of matrices with entries in a ring $A$. The Cohn localization of $A$ with respect to $W$ is a ring obtained from $A$ by formally inverting all
matrices in $W$. This ring $\Lambda$ is equipped with a ring homomorphism $A\rightarrow\Lambda$ and every matrix in $W$ becomes invertible in $\Lambda$. Moreover
$\Lambda$ is universal with respect to this property.

A particular example of such a localization may be obtained from a ring homomorphism $A\rightarrow B$. In this case the Cohn localization of $A$ with respect to the set
of matrices with entries in $A$ becoming invertible in $B$ will be denoted by $L(A\rightarrow B)$. Because of the universal property, the morphism $A\rightarrow B$
factors through $L(A\rightarrow B)$.

If $A$ is commutative, a Cohn localization $\Lambda$ of $A$ is a classical localization $S^{-1}A$, where $S$ is a multiplicative set in $A$ and $\Lambda$ is a flat
$A$-module. In the general case $\Lambda$ is much more complicated, in
particular if $A\rightarrow B$ is the augmentation map of a group ring $\Z[G]$. In this case, the ring $\Lambda=L(\Z[G]\rightarrow\Z)$ is known if $G$ is commutative or
free [FV]. See also [CR] for an explicit description of $\Lambda$ when $G$ is cyclic. But in the other cases almost nothing is known.

The Cohn localization is an interesting construction in algebra, but this theory plays an important role in homotopy theory, algebraic K-theory and homology surgery.

The main example of a geometric problem involving the Cohn localization is the problem of codimension $2$ embeddings. If we want to classify the embeddings from
a closed maniflod $M$ to a closed manifold $W$, we need to determine the complement $X$ of this embedding. The homology of $X$ with coefficients in $\Z[\pi_1(W)]$ can be
determined by the geometry of the embedding. If the codimension of the embedding is more than $2$, the fundamental group of $X$ is isomorphic to $\pi_1(W)$ and the
homotopy type of $X$ can be determined (in some sense). But, in the codimension $2$ case, the situation is completely different because the morphism
$\pi_1(X)\rightarrow\pi_1(W)$ is never an isomorphism.

There is a strategy to study the problem in this case (see [CS] and [LD]):

First of all, we have to choose a group $\pi$ and a morphism $\pi\rightarrow\pi_1(W)$ which is possibly a morphism $\pi_1(X)\rightarrow\pi_1(W)$ obtained from an
embedding. Then we set: $A=\Z[\pi]$, $B=\Z[\pi_1(W)]$. We have a ring homomorphism $A\rightarrow B$ and then a localization ring $\Lambda=L(A\rightarrow B)$.

Second step: for a possible space $X$ and a possible map $X\rightarrow W$, we have to find a space $Y$ between $X$ and $W$ which is universal with respect to the
following property:

$\bullet$ the map $X\rightarrow Y$ induces isomorphisms $\pi_1(X)\build\longrightarrow_{}^\sim\pi_1(Y)$ and $\H_*(X,B)\build\longrightarrow_{}^\sim\H_*(Y,B)$.

Such a space $Y$ is obtained by localization in the category of spaces over $W$ and the homotopy groups of $Y$ are $\Lambda$-modules.

Last step: find a manifold $M'$ and a $B$-homology equivalence $M'\rightarrow Y$. The obstruction to finding such a manifold lies, roughly speaking, in a surgery group 
$\Gamma_n(A\rightarrow B)$ defined by Cappell and Shaneson in [CS] which is isomorphic to the classical surgery obstruction group $L_n(\Lambda)$ [V2].

If all this procedure works $M'$ will be the complement of a tubular neighbourhood of an embedding $M\subset W$.

The Cohn localization ring $\Lambda$ appears twice in this theory: in the homotopy groups of $Y$ and in the obstruction surgery group.
\vskip 12pt
The first result of this paper is a complete description of $L(\Z[G]\rightarrow\Z)$ when $G$ is a finite group.
\vskip 12pt
\noi{\bf Theorem A.} {\sl Let $G$ be a finite group. For every prime $p$, denote by $G_p$ the quotient of $G$ by the subgroup of $G$ generated by all the elements
of order coprime to $p$, and by $\varepsilon_p:\Z_{(p)}[G_p]\rightarrow\Z_{(p)}$ the corresponding augmentation map. Then the Cohn localization ring
$\Lambda=L(\Z[G]\rightarrow\Z)$ is given by the following pull-back diagram:
\vskip -24pt
$$\diagram{\Lambda&\hfl{}{}&\build\prod_p^{}\Z_{(p)}[G_p]\cr \vfl{}{}&&\vfl{}{\prod_p \varepsilon_p}\cr \Z&\hfl{\Delta}{}&\build\prod_p^{}\Z_{(p)}\cr}$$
\vskip -24pt 
\noi where $\Delta$ is the diagonal inclusion and the product is over all non-trivial $G_p$.}
\vskip 12pt
Consider a ring homomorphism $A\rightarrow B$. We'll say that $B$ is a central localization (resp. an Ore localization) of $A$ if $B$ is the ring $S^{-1}A$ where $S$ is
a multiplicative set in the center of $A$ (resp. a multiplicative set in $A$ satisfying the Ore condition). 

We say that $B$ is stably flat over $A$ if the two conditions hold:

\noi$\bullet$ the multiplication map: $B\build\otimes_A^{}B\rightarrow B$ is an isomorphism

\noi$\bullet$ Tor$^A_i(B,B)=0$ for all $i>0$.
\vskip 12pt
The second result is the following:
\vskip 12pt
\noi{\bf Theorem B.} {\sl Let $G$ be a finite group and $\Lambda$ be the  Cohn localization ring $L(\Z[G]\rightarrow\Z)$. Then the following conditions are equivalent:

1) $\Lambda$ is a central localization of $\Z[G]$

2) $\Lambda$ is an Ore localization of $\Z[G]$

3) $\Lambda$ is a flat left $\Z[G]$-module

4) $\Lambda$ is stably flat over $\Z[G]$

5) $G$ is nilpotent.}
\vskip 12pt
Notice that 2) and 3) are equivalent in any case by a result of Teichner [T].
\vskip 12pt
The section 1 is devoted to a brief presentation of the Cohn localization and some useful properties of this functor.

In section 2, the localization functor is extended to the category of (chain) complexes over a ring $A$. The universal property of this localization functor $\Phi$ and
the ring structure of $A$ induce a graded ring structure on the homology $\Lambda_*$ of the localization complex $\Phi(A)$.
The subring $\Lambda_0$ of the graded localization ring $\Lambda_*$ is actually the Cohn localization of $A$. Moreover, for each complex $C$, there is a spectral
sequence converging to $\H_*(\Phi(C))$, where the $E^2$ term is given by: $E^2_{pq}=\H_p(C\otimes_A\Lambda_q)$.

The last section is devoted to the Cohn localization ring $\Lambda=L(\Z[G]\rightarrow\Z)$ when $G$ is any finite group. Using results of sections 1 and 2, Theorem A and
Theorem B are proven and a complete description of the graded localization ring $\Lambda_*$ is given (Theorem 3.8). There is also some examples and counter-examples at
the end of the paper.
\vskip 24pt
\noi{\bf 1. Cohn localization of rings.}
\vskip 12pt
For any ring $A$, the class of finitely generated projective right $A$-modules will be denoted by ${\rde P}(A)$. The set of matrices with entries in $A$ will be denoted
by $M(A)$ and the class of morphisms between two modules in ${\rde P}(A)$ will be denoted by ${\rde M}(A)$. The set $M(A)$ may be seen as a subset of ${\rde M}(A)$.

Let $A$ be a ring and $W$ be a class of morphisms in ${\rde M}(A)$. The ring obtained from $A$ by formally inverting all morphisms in $W$ is called the Cohn localization
(or the universal localization) of $A$ with respect to $W$. This localization ring $\Lambda$ is equipped with a ring homomorphism from $A$ to
$\Lambda$ and the morphism $A\rightarrow\Lambda$ is characterized by the following properties:

$\bullet$ every morphism in $W$ is sent to an isomorphism in ${\rde M}(\Lambda)$.

$\bullet$ $\Lambda$ is universal with respect to this property.

This localization ring will be denoted by $L(A,W)$ or $W^{-1}A$.
\vskip 12pt
This universal characterization implies that the ring $L(A,W)$ is unique up to canonical isomorphism if it exists. A construction of $L(A,W)$ by generators and relations
is given in [C1], [C2] if $W$ is a set or matrices. Actually this construction may be modified in order to produce $L(A,W)$ in the general case. See also a presentation
of the Cohn localization in [R].

On the other hand, a construction of this ring will be obtained as a consequence of Theorem 2.7 in the next section.
\vskip 12pt
Let $f:A\rightarrow B$ be a ring homomorphism. There is an exact sequence in algebraic K-theory:
$$\cdots\longrightarrow K_1(A)\longrightarrow K_1(B)\longrightarrow K_1(f)\build\longrightarrow_{}^\partial K_0(A)\longrightarrow K_0(B)\longrightarrow \cdots$$
and the relative K-theory group $K_1(f)=K_1(B,A)$ is generated by triples $(P,Q,\lambda)$, where $P$ and $Q$ are in ${\rde P}(A)$ and $\lambda$ is an isomorphism from
$P\otimes_A B$ to $Q\otimes_A B$.

If $(P,Q,\lambda)$ is such a triple, the class of $(P,Q,\lambda)$ in $K_1(f)$ will be denoted by $\tau(P,Q,\lambda)$. 

Let $u$ be a subgroup of $K_1(f)$. The class of morphisms $\alpha:P\rightarrow Q$ in ${\rde M}(A)$ such that $\alpha\otimes_A B$ is an isomorphism and
$\tau(P,Q,\alpha\otimes_A B)$ belongs to $u$ will be denoted by $W^u(A\rightarrow B)$ and the localization ring of $A$ with respect to this class will be denoted by
$L^u(A\rightarrow B)$.

If $u$ is the trivial subgroup of $K_1(f)$ (resp. the all group $K_1(f)$, the image of $K_1(B)\rightarrow K_1(f)$), the localization ring $L^u(A\rightarrow B)$ will be
denoted by $L^s(A\rightarrow B)$ (resp. $L^p(A\rightarrow B)$, $L^h(A\rightarrow B)$).

If $u\subset v$ are two subgroups of $K_1(f)$, we have: $W^u(A\rightarrow B)\subset W^v(A\rightarrow B)$ and there is a canonical ring homomorphism:
$L^u(A\rightarrow B)\rightarrow L^v(A\rightarrow B)$.

If the ring $L^u(A\rightarrow B)$ doesn't depend of $u$ it will be denoted by $L(A\rightarrow B)$.
\vskip 12pt
\noi{\bf Remark:} Actually every Cohn localization ring $L(A,W)$ is the localization ring $L^u(A\rightarrow B)$ for some ring homomorphism $A\rightarrow B$ and some
subgroup of $K_1(A\rightarrow B)$. More precisely we have the following:

Let $A$ be a ring and $W$ be a set of morphisms in ${\rde M}(A)$. Let $\Lambda$ be the localization ring $L(A,W)$ and $u$ be a subgroup
of $K_1(A\rightarrow\Lambda)$ containing all the elements $\tau(P,Q,\alpha\otimes_A\Lambda)$ with $\alpha:P\rightarrow Q$ in $W$. Then $\Lambda$ is the localization ring
$L^u(A\rightarrow\Lambda)$.
\vskip 12pt
\noi{\bf 1.1 Lemma:} {\sl Let $f:A\rightarrow B$ be a ring homomorphism. Let $W^h$ be the set of square matrices in $M(A)$ becoming invertible in $M(B)$ and
$W^s$ be the set of square matrices $M\in M(A)$ such that $f(M)$ is an invertible matrix with zero torsion in $K_1(B)$. Then we have: $L^h(A\rightarrow B)=L(A,W^h)$
and $L^s(A\rightarrow B)=L(A,W^s)$.}
\vskip 12pt
\noi{\bf Proof:} Denote by $\Lambda^s$ and $\Lambda^h$ the localization rings $L(A,W^s)$ and $L(A,W^h)$.

Let $u$ be the trivial group in $K_1(f)$ and $v$ be the image of $K_1(B)\rightarrow K_1(f)$. The sets $W^h$ and $W^s$ are included in
$W^u(A\rightarrow B)$ and $W^v(A\rightarrow B)$. Then we have ring homomorphisms $\Lambda^s\rightarrow L^s(A\rightarrow B)$ and $\Lambda^h\rightarrow L^h(A\rightarrow B)$
and the only thing to do is to prove that every morphism in $W^u(A\rightarrow B)$ (resp. $W^v(A\rightarrow B)$) is sent to an isomorphism in ${\rde M}(\Lambda^s)$ (resp.
in ${\rde M}(\Lambda^h)$).

Let $\alpha:P\rightarrow Q$ be a morphism in $W^v(A\rightarrow B)$ and $\theta$ be the class $\tau(P,Q,\alpha\otimes B)$ in $K_1(f)$. Then $\alpha$ belongs
to  $W^u(A\rightarrow B)$ if and only if $\theta=0$. Since $\theta$ belongs to the image of $K_1(B)\rightarrow K_1(f)$, we have: $\partial(\theta)=0$ in $K_0(A)$ and
$P$ and $Q$ are equal in $K_0(A)$. Therefore there exists a module $K\in {\rde P}(A)$ such that $P\oplus K$ and $Q\oplus K$ are isomorphic and $P\oplus K$ is free. Let
$\beta$ be the morphism $\alpha\oplus\hbox{Id}_K$. Then $\beta$ is represented by a square matrix in $W^h$ and $\beta$ is invertible in $M(\Lambda^h)$. Hence $\alpha$ is
sent to an isomorphism in ${\rde M}(\Lambda^h)$.

Suppose $\theta=0$. Then $\alpha$ belongs to $W^u(A\rightarrow B)$ and $\beta$ also. There exist a free module $F\in{\rde P}(A)$ and two isomorphisms
$\varepsilon: F\build\rightarrow_{}^\sim P\oplus K$ and $\varepsilon':F\build\rightarrow_{}^\sim Q\oplus K$ and we have:
$$\tau(P\oplus K,Q\oplus K,\beta\otimes B)=\tau(F,F,(\varepsilon'^{-1}\circ\beta\circ\varepsilon)\otimes B)=0$$
Set: $\beta'=\varepsilon'^{-1}\circ\beta\circ\varepsilon$. Then the torsion $\tau'$ of $\beta'\otimes B$ in $K_1(B)$ vanishes in $K_1(f)$ and $\tau'$ comes
from an element $\widehat\tau\in K_1(A)$. Up to stabilize $K$ and $F$ we may as well suppose that $\widehat\tau$ is the torsion of an automorphism $\widehat\varepsilon$
in GL$(F)$. Consider the morphism: $\widehat\beta=\widehat\varepsilon^{-1}\circ\beta'$. This morphism corresponds to a square matrix with zero torsion and
$\widehat\beta$ is invertible in $M(\Lambda^s)$. Therefore $\beta$ and $\alpha$ are sent to isomorphisms in ${\rde M}(\Lambda^s)$.\cqfd 
\vskip 12pt
\noi{\bf 1.2 Examples:} Let $A$ be a commutative ring and $W$ be a set of matrices in $M(A)$. Then we have the following:

If $W$ contains some non-square matrix then the ring $L(A,W)$ is the trivial ring.

If $W$ contains only square matrices the ring $L(A,W)$ is the localized ring $S^{-1}A$, where $S$ is the multiplicative set in $A$ generated by the determinants of the
matrices in $W$.

Let $f:A\rightarrow B$ be a ring homomorphism between commutative rings. Consider the following subsets of $A$: the set $S=f^{-1}(1)$ and the set $\Sigma=f^{-1}(B^*)$,
where $B^*$ is the group of units in $B$. Then we have:
$$L^s(A\rightarrow B)=S^{-1} A$$
$$L^h(A\rightarrow B)=\Sigma^{-1} A$$

\vskip 12pt
In the non commutative case, the localization ring is much more difficult to understand. Nevertheless it is possible to give some description of morphisms in
${\rde M}(\Lambda)$. In order to do that, we need some notations:
\vskip 12pt
\noi{\bf Notation}: Let $A$ be a ring and $\alpha$, $\beta_1$, $\beta_2$,\dots,$\beta_p$ be some morphisms in ${\rde M}(A)$. We say that $\alpha$ is a compatible product
$\beta_1\beta_2\dots\beta_n$ if there are finitely generated projective right $A$-modules  $P_0,P_1,\dots,P_n$ such that $\beta_i$ is a morphism from $P_i$ to $P_{i-1}$,
$\alpha$ is a morphism from $P_n$ to $P_0$ and $\alpha$ is equal to the composite $\beta_1\beta_2\dots\beta_p$.
\vskip 12pt
\noi{\bf Notation:} {\sl Let $A$ be a ring and $W$ be a class in ${\rde M}(A)$. The class of morphisms in ${\rde M}(A)$ of the form:
$$\pmatrix{\alpha_1&*&.&.&.&*\cr 0&\alpha_2&.&.&.&*\cr .&.&.&&&.\cr .&.&&.&&.\cr .&.&&&.&*\cr 0&0&.&.&0&\alpha_p\cr}$$
where each $\alpha_i$ is in $W$ and each $*$ is in ${\rde M}(A)$, will be denoted by $\widehat W$. 

The classes $W$ and $\widehat W$ induce the same localization ring.}
\vskip 12pt
\noi{\bf 1.3 Proposition:} {\sl Let $A$ be a ring, $W$ be a class of morphisms in ${\rde M}(A)$ and $\Lambda$ be the ring $L(A,W)$. Denote by $M_0$ the image of
${\rde M}(A)$ in ${\rde M}(\Lambda)$ and by $M_1$ the image of $\widehat W$ in ${\rde M}(\Lambda)$. Every morphism in $M_1$ is clearly invertible.

Then for every $P$ and $Q$ in ${\rde P}(A)$, every morphism from $P\otimes_A\Lambda$ to $Q\otimes_A\Lambda$ is a compatible product
$\lambda\alpha^{-1}\mu$ with $\lambda$ and $\mu$ in $M_0$ and $\alpha$ in $M_1$.} 
\vskip 12pt
\noi{\bf Proof:} Let ${\rde C}$ be the class of all $\Lambda$-modules on the form $M\otimes_A\Lambda$ with $M\in{\rde P}(A)$
and ${\rde M}$ be the class of morphisms in ${\rde M}(\Lambda)$ between two modules in ${\rde C}$.
Denote by ${\rde M}'(\Lambda)$ the class of compatible products $\lambda\alpha^{-1}\mu$ with $\lambda,\mu$ in $M_0$ and $\alpha$ in $M_1$. Clearly, ${\rde M}'(\Lambda)$
is contained in ${\rde M}$. Because of the following:
$$\lambda\alpha^{-1}\mu+\lambda'\alpha'^{-1}\mu'=\pmatrix{\lambda&\lambda'\cr}\pmatrix{\alpha&0\cr0&\alpha'\cr}^{-1}\pmatrix{\mu\cr \mu'\cr}$$
the class ${\rde M}'(\Lambda)$ is stable under (compatible) sum. Moreover we have:
$$(\lambda\alpha^{-1}\mu)(\lambda'\alpha'^{-1}\mu')=\pmatrix{-\lambda&0\cr}\pmatrix{\alpha&\mu\lambda'\cr0&\alpha'\cr}^{-1}\pmatrix{0\cr\mu'\cr}$$
and ${\rde M}'(\Lambda)$ is stable under compatible product.

In particular the set of $1\times1$-matrices in ${\rde M}'(\Lambda)$ is a subring $\Lambda_0$ of $\Lambda$ containing the image of $A\rightarrow\Lambda$.
Therefore every matrix in ${\rde M}'(\Lambda)$ is a matrix in ${\rde M}(\Lambda_0)$ and ${\rde M}'(\Lambda)$ is contained in ${\rde M}(\Lambda_0)$.

On the other hand every morphism in the image of $W\rightarrow{\rde M}'(\Lambda)$ is invertible in ${\rde M}'(\Lambda)\subset{\rde M}(\Lambda_0)$ and the ring $\Lambda_0$
satisfies the universal property of $\Lambda$. Hence one has: $\Lambda=\Lambda_0$ and the class ${\rde M}'(\Lambda)$ is the class of morphisms in ${\rde M}(\Lambda)$
between two $\Lambda$-modules in ${\rde C}$. The result follows.\cqfd
\vskip 12pt
\noi{\bf 1.4 Corollary:} {\sl Let $A\rightarrow B$ be a ring homomorphism and $u$ be a subgroup of $K_1(A\rightarrow B)$ . Let $\Lambda$ be the Cohn localization ring
$L^u(A\rightarrow B)$. Then for every module $P\in{\rde P}(\Lambda)$ and every morphism $\varphi: P\rightarrow P$, the following holds:
$$\varphi\otimes_\Lambda B=\hbox{Id}_{P\otimes B}\ \ \Longrightarrow\ \ \ \varphi\ \hbox{is an isomorphism.}$$}
\vskip 12pt
\noi{\bf Proof:} Let $W$ be the class $W^u(A\rightarrow B)$. The class $\widehat W$ is actually the class $W$.

Let $\varphi:P\rightarrow P$ be a morphism in ${\rde M}(\Lambda)$ such that $\varphi\otimes_\Lambda B$ is the identity. We have to prove that $\varphi$ is an isomorphism.

First of all, suppose $P$ is free. Then there is a free module $F\in{\rde P}(A)$ such that: $P=F\otimes_A\Lambda$ and $\varphi\otimes_\Lambda B$ is the identity of
$F\otimes_A B$.

Because of Proposition 1.3, there exist two modules $P_1$ and $Q_1$ in ${\rde P}(A)$ and three morphisms:
$$F\build\longrightarrow_{}^\lambda P_1\build\longleftarrow_{}^\alpha Q_1\build\longrightarrow_{}^\mu F$$
such that $\alpha$ belongs to $W$ and we have: $\varphi=\lambda\otimes_A\Lambda\circ(\alpha\otimes_A\Lambda)^{-1}\circ\mu\otimes_A\Lambda$.

Denote by $\lambda'$, $\alpha'$, $\mu'$ (resp. $\lambda''$, $\alpha''$, $\mu''$) the images of $\lambda$, $\alpha$, $\mu$ by the tensorization functor $\otimes_A\Lambda$
(resp. $\otimes_A B$).  Then we have: $\varphi=\lambda'\alpha'^{-1}\mu'$.

Let $\psi$ be the morphism $\pmatrix{0&\lambda\cr\mu&\alpha\cr}$ from $F\oplus Q_1$ to $F\oplus P_1$. By tensorization with $\Lambda$ and $B$, we get the morphisms
$\psi'$ and $\psi''$.

Since $\varphi\otimes_A B=\lambda''\alpha''^{-1}\mu''$ is the identity, the morphism $\psi''$ is invertible with inverse:
$$\psi''^{-1}=\pmatrix{0&\lambda''\cr\mu''&\alpha''\cr}^{-1}=\pmatrix{-1&\lambda''\alpha''^{-1}\cr \alpha''^{-1}\mu''&(1-\alpha''^{-1}\mu''\lambda'')\alpha''^{-1}\cr}$$

An elementary computation shows the following in $K_1(A\rightarrow B)$:
$$\tau(F\oplus Q_1,F\oplus P_1,\psi'')=\tau(Q_1,P_1,\alpha'')$$
and $\tau(F\oplus Q_1,F\oplus P_1,\psi'')$ belongs to $u$. By the universal property of $\Lambda$, $\psi'$ is invertible.

Let $\pmatrix{a&b\cr c&d\cr}$ be the inverse of $\psi'$. We have the following relations:
$$\lambda' c=1\hskip 24pt  \lambda' d=0\hskip 24pt  \mu' a+\alpha' c=0\hskip 24pt  \mu' b+\alpha' d=1$$
$$b\mu'=1\hskip 24pt  a\lambda'+b\alpha'=0\hskip 24pt  d\mu'=0\hskip 24pt  c\lambda'+d\alpha'=1$$
So we have:
$$0=\lambda'\alpha'^{-1}(\mu' a+\alpha' c)=\varphi a+\lambda' c=\varphi a+1$$
$$0=(a\lambda'+b\alpha')\alpha'^{-1}\mu'=a\varphi+b\mu'=a\varphi+1$$
Therefore $\varphi$ is invertible (with inverse $-a$).

Consider now the general case. 

Let $Q$ be a module in ${\rde P}(\Lambda)$ such that $P\oplus Q$ is free. Because of the previous study, $\varphi\oplus\hbox{Id}_Q$ is invertible and $\varphi$
is invertible too. The result follows.\cqfd 
\vskip 12pt
\noi{\bf 1.5 Proposition:} {\sl Let $A\rightarrow B$ be a ring homomorphism and $u$ be a subgroup of $K_1(A\rightarrow B)$. Let 
$\Lambda=L^u(A\rightarrow B)$ be the corresponding Cohn localization ring. Suppose the induced morphism
$\Lambda\rightarrow B$ is surjective. Then, for every finitely generated right $\Lambda$-module $M$, the following holds:
$$M\otimes_\Lambda B=0\ \ \Longrightarrow\ \ M=0$$}
\vskip 12pt
\noi{\bf Proof:} Suppose $M$ is a finitely generated right $\Lambda$-module and $M\otimes_\Lambda B$ is trivial. Take a free resolution of $M$:
$$C_1\build\longrightarrow_{}^d C_0\longrightarrow M\longrightarrow 0$$
where $C_0$ is finitely generated. Since $M\otimes_\Lambda B=0$ the morphism 
$$C_1\otimes_\Lambda B\build\longrightarrow_{}^d C_0\otimes_\Lambda B$$
is surjective and has a section $s$. Since $\Lambda\rightarrow B$ is surjective $s$ can be lifted to a morphism $\widetilde s: C_0\rightarrow C_1$.
The composite morphism $d\circ\widetilde s$ corresponds to a matrix in $M(\Lambda)$ which is sent to an identity matrix in $M(B)$. Because of Corollary 1.4, $u$ is
invertible and the map $d: C_1\longrightarrow C_0$ is surjective. The result follows.\cqfd
\vskip 12pt
\noi{\bf 1.6 Corollary:} {\sl Let $A\rightarrow B$ be a ring homomorphism and $u$ be a subgroup of $K_1(A\rightarrow B)$. 

Suppose the morphism $L^u(A\rightarrow B)\rightarrow B$ is onto. Then, for every subgroup $v$ of $K_1(A\rightarrow B)$ containing $u$ the canonical
ring homomorphism from $L^u(A\rightarrow B)$ to $L^v(A\rightarrow B)$ is an isomorphism.}
\vskip 12pt
\noi{\bf Proof:} Suppose $L^u(A\rightarrow B)\longrightarrow B$ is onto. Let $\varphi:P\rightarrow Q$ be a morphism in ${\rde M}(\Lambda)$ such that
$\varphi\otimes_\Lambda B$ is an isomorphism. Let $K$ and $C$ be the kernel and the cokernel of $\varphi$. The module $C$ if finitely generated and $C\otimes_\Lambda B$
is trivial. Because of Proposition 1.5, $C$ is trivial too and $\varphi$ is surjective. Hence the sequence:
$$0\longrightarrow K\longrightarrow P\longrightarrow Q\longrightarrow0$$
is split exact and $K$ is finitely generated. Then we have a split exact sequence:
$$0\longrightarrow K\otimes_\Lambda^{}B\longrightarrow P\otimes_\Lambda^{}B\build\longrightarrow_{}^\sim Q\otimes_\Lambda^{}B\longrightarrow0$$
and the module $K\otimes_\Lambda B$ is trivial. Because of Proposition 1.5 the module $K$ is also trivial and $\varphi$ is an isomorphism.

Therefore $L^u(A\rightarrow B)$ satisfies the universal property of $L^v(A\rightarrow B)$ and the result follows.\cqfd
\vskip 12pt
\noi{\bf 1.7 Remark:} If the morphism $A\rightarrow B$ is onto, the ring $L^u(A\rightarrow B)$ doesn't depend on $u$ and will be denoted by $L(A\rightarrow B)$.
\vskip 12pt
\noi{\bf 1.8 Corollary:} {\sl Let $G$ be a group and $H$ be a finitely generated perfect subgroup of $G$. Let $N$ be the normalizer of $H$. Then the morphism
$\Z[G]\rightarrow L(\Z[G]\rightarrow\Z)$ factors through $\Z[G/N]$ and $L(\Z[G]\rightarrow\Z)$ is the ring $L(\Z[G/N]\rightarrow\Z)$.}
\vskip 12pt
\noi{\bf Proof:} Denote by $\Lambda$ the localization ring $L(\Z[G]\rightarrow\Z)$ and by $f:\Z[G]\rightarrow\Lambda$ the induced morphism. Let $I$ be the kernel of the
morphism $\Z[G]\longrightarrow \Z[H\backslash G]$. Since $H$ is finitely generated, $I$ is a finitely generated right $\Z[G]$-module and $J=I\otimes_{\Z[G]}\Lambda$ is a
finitely generated right $\Lambda$-module. We have an exact sequence:
$$0\longrightarrow I\longrightarrow \Z[G]\longrightarrow\Z[H\backslash G]\rightarrow 0$$
and that implies:
$$I\otimes_{\Z[G]}\Z\simeq \hbox{Tor}_1^{\Z[G]}(\Z[H\backslash G],\Z)\simeq \H_1(H)\simeq 0$$
$$\Longrightarrow\ \ J\otimes_\Lambda\Z=I\otimes_{\Z[G]}\Lambda\otimes_\Lambda\Z=I\otimes_{\Z[G]}\Z=0$$
Hence: $J=0$ because of Proposition 1.5. 

On the other hand, the tensorization by $\Lambda$ induces a commutative diagram:
$$\diagram{I&\hfl{}{}&\Z[G]\cr\vfl{}{}&&\vfl{}{}\cr J&\hfl{}{}&\Lambda\cr}$$
Since $J$ is trivial, $I$ vanishes in $\Lambda$. Thus, for every $h$ in $H$, $h-1$ vanishes in $\Lambda$ and we have: $f(H)=\{1\}$ and therefore: $f(N)=\{1\}$.
Hence the morphism $f$ factors through $\Z[G/N]$. Moreover the morphism $\Z[G/N]\rightarrow\Lambda$ satisfies the universal condition defining the ring 
$L(\Z[G/N]\rightarrow\Z)$ and $\Lambda$ is the localization ring $L(\Z[G/N]\rightarrow\Z)$.\cqfd
\vskip 12pt
\noi{\bf 1.9 Lemma:} {\sl Let $f:A\rightarrow B$ be a surjective ring homomorphism. Then the map $A\rightarrow L(A\rightarrow B)$ is an isomorphism if and only if 
every element in $f^{-1}(1)$ is invertible.}
\vskip 12pt
\noi{\bf Proof:} Denote by $\Lambda$ the localization ring $L(A\rightarrow B)$. Suppose $A\rightarrow\Lambda$ is an isomorphism. Let $a$ be an element in $f^{-1}(1)$.
This element represents a $1\times 1$-matrix which is invertible in $M(B)$ and then in $M(\Lambda)$. Because $A\rightarrow\Lambda$ is bijective, the matrix is
invertible. Thus $a$ is invertible and $f^{-1}(1)$ is contained in $A^*$ (the group of units of $A$).

Suppose now $f^{-1}(1)$ is contained in $A^*$. Consider a square matrix $M$ in $M(A)$ such that $f(M)$ is invertible. Let $N$ be a matrix in $M(A)$ with $f(N)=f(M)^{-1}$.
Since the diagonal entries of $MN$ are invertible in $A$ and the other entries are killed in $B$, it is possible to multiply $MN$ on the left and the right by
elementary matrices in order to obtain a diagonal matrix
and the same holds for $NM$. Therefore $MN$ and $NM$ are invertible and $M$ is invertible too. Then every square matrix in $M(A)$ sent to an invertible matrix in $M(B)$
is invertible and, because of Lemma 1.1, we have:
$$L(A\rightarrow B)=L^h(A\rightarrow B)=A\eqno{\cqfd}$$
\vskip 12pt
As a corollary one gets the following result:
\vskip 12pt
\noi{\bf 1.10 Proposition:} {\sl Let $p$ be a prime and $G$ be a finite $p$-group. Denote by $\Lambda$ the localization ring $L(\Z_{(p)}[G]\rightarrow\Z_{(p)})$. Then the
map $\Z_{(p)}[G]\rightarrow\Lambda$ is an isomorphism.}
\vskip 12pt
\noi{\bf Proof:} Let's say that a ring homomorphism $f$ is local if $f$ is surjective and every element in $f^{-1}(1)$ is invertible. 

The proposition is obvious if the order of $G$ is $1$. The result will be proven by induction. Suppose the order of $G$ is $p^n>1$ and the proposition is
true for every group of order $p^i$, $i<n$. Let $z\in G$ be a central element of order $p$ and $G'$ be the quotient $G/<z>$. Consider the maps:
$$\Z_{(p)}[G]\build\longrightarrow_{}^f \Z_{(p)}[G']\build\longrightarrow_{}^g\Z_{(p)}$$
Because of lemma 1.9 $g$ is local and we have to prove that $g\circ f$ is also local. 

Let $f'$ be the reduction of $f$ mod $p$:
$$f':\F_p[G]\longrightarrow \F_p[G']$$
Set: $z=1+u$. The relation: $z^p=1$ becomes mod $p$: $u^p=0$. Let $U$ be an element in $f'^{-1}(1)$. This element has the form: $U=1-uV$ for some element $V\in \F_p[G]$.
Then $U$ is invertible with inverse:
$$U^{-1}=1+uV+u^2V^2+\dots+u^{p-1}V^{p-1}$$
and $f'$ is local. Let $U$ be an element of $f^{-1}(1)\subset\Z_{(p)}[G]$. The multiplication by $U$ is an endomorphism $\varphi$ on $\Z_{(p)}[G]$ considered as a free
$\Z_{(p)}$-module. Because $U$ is invertible mod $p$, the determinant of $\varphi$ is non zero mod $p$ and then invertible in $\Z_{(p)}$. Consequently, $U$ is invertible
and $f$ and $g\circ f$ are local. The result follows.\cqfd
\vskip 12pt
We have a last result in this section which will be used in the sequel.
\vskip 12pt
\noi{\bf 1.11 Proposition:} {\sl Let $f:A\rightarrow B$ and $g:B\rightarrow R$ be two ring homomorphisms. Suppose $f$ is onto. Let $I$ be the kernel of $f$ and
$\Lambda$ be the Cohn localization $L^h(A\rightarrow R)$. Let $J$ be the two-sided ideal of $\Lambda$ generated by the image of $I$ in $\Lambda$. Then the Cohn
localization $L^h(B\rightarrow R)$ is the quotient $\Lambda/J$.}
\vskip 12pt
\noi{\bf Proof:} Since $I\subset A$ is killed in $R$, $J$ is also killed in $R$ and we have a commutative diagram of rings:
$$\diagram{A&\hfl{f}{}&B&\hfl{g}{}&R\cr \vfl{}{}&&\vfl{}{}&&\vfl{}{\hbox{Id}}\cr\Lambda&\hfl{}{}&\Lambda/J&\hfl{}{}&R\cr}$$
Let $M$ be a square matrix in $M(B)$ which is invertible in $M(R)$. Since $f$ is onto, $M$ can be lifted in a matrix $N\in M(A)$ which is invertible in $M(\Lambda)$ and
then in $M(\Lambda/J)$. Therefore $M$ is invertible in $M(\Lambda/J)$. 

Let $h:B\rightarrow C$ be a ring homomorphism such that every square matrix in $M(B)$ becoming invertible in $M(R)$ is invertible in $M(C)$. Be the universal property,
the morphism $h\circ f:A\rightarrow C$ factors through $\Lambda$. Moreover $J$ is killed by the induced morphism $\Lambda\rightarrow C$ and this morphism factors
through $\Lambda/J$. So we have the following commutative diagram of rings:
$$\diagram{A&\hfl{f}{}&B&\hfl{g}{}&C\cr \vfl{}{}&&\vfl{}{}&&\vfl{}{\hbox{Id}}\cr\Lambda&\hfl{}{}&\Lambda/J&\hfl{}{}&C\cr}$$
and the morphism $B\rightarrow C$ factors in a unique way through $\Lambda/J$. By the universal property, $\Lambda/J$ is the Cohn localization $L^h(B\rightarrow R)$.\cqfd
\vskip 24pt
\noi{\bf 2. Localization of complexes.}
\vskip 12pt
This section is devoted to the localization in the category of complexes. For similar theories see also [V2] section 5, [D], [R], [NR].
\vskip 12pt 
Let $R$ be a ring. The graded differential right $R$-modules which are projective and bounded from below define a category ${\rde C}_*(R)$. The objects in
${\rde C}_*(R)$ are called complexes (or $R$-complexes if needed). A morphism in ${\rde C}_*(R)$ is a linear map of some degree which commutes with the differentials
(in the graded sense). A morphism in ${\rde C}_*(R)$ is a cofibration if it is injective with cokernel in ${\rde C}_*(R)$. A complex $C$ is said to be finite if it is
bounded, finitely generated and free in each degree. A complex $C$ is said to be finitely dominated (or compact) if $C$ is bounded and finitely generated.

A complex $C$ is finitely dominated if and only if it is a direct summand of a finite complex.

Let $C$ and $C'$ be two $R$-complexes. Denote by $\Hom(C,C')_p$ the module of $R$-linear maps of degree $p$ from $C$ to $C'$. The direct sum of all these modules is a
($\Z$--)graded module denoted by $\Hom(C,C')$. There is a unique differential $d$ on this graded module such that:
$$\forall (x,f)\in C\times\Hom(C,C'),\ \ \ d(f(x))=d(f)(x)+(-1)^{|f|} f(d(x))$$
where $|f|$ is the degree of $f$. Equipped with this differential, $\Hom(C,C')$ is a graded differential module. A cycle in $\Hom(C,C')$ is a morphism and a boundary
is a homotopy. Moreover if $C$, $C'$ and $C''$ are $R$-complexes, the composition map:
$$\Hom(C',C'')\otimes\Hom(C,C')\longrightarrow\Hom(C,C'')$$
is a degree $0$ morphism of graded differential modules.
\vskip 12pt
\noi{\bf 2.0 Lemma:} {\sl Let $C$ and $C'$ be two $R$-complexes. Then one has the following properties:

$\bullet$ the graded differential module $\Hom(C,C')$ is acyclic if and only if every morphism from $C$ to $C'$ is null-homotopic.

$\bullet$ the graded differential module $\Hom(C,C)$ is a graded differential ring and $\Hom(C,C')$ is a graded differential right $\Hom(C,C)$-module.}
\vskip 12pt
\noi{\bf Proof:} It's just a straightforward computation.\cqfd
\vskip 12pt

Throughout this section $A$ is a ring and $W$ is a class of morphisms in ${\rde M}(A)$. The Cohn localization $L(A,W)$ will be denoted by $\Lambda$.

The class of complexes of length $2$:
$$\dots\longrightarrow 0\longrightarrow P\build\longrightarrow_{}^d Q\longrightarrow 0\longrightarrow\dots$$
such that the only non zero differential $d$ lies in $W$ will be denoted by ${\rde W}_0$.
\vskip 12pt
\noi{\bf Definitions: } {\sl A right $A$-module $M$ is said to be $W$-local if, for every morphism $f: P\rightarrow Q$ in $W$,
the induced map $f^*: \Hom(Q,M)\rightarrow\Hom(P,M)$ is bijective.

An $A$-complex $C$ is said to be $W$-local if every morphism from a complex in ${\rde W}_0$ to $C$ is null-homotopic.

Let ${\rde C}$ be a class of $A$-complexes. An $A$-complex $C$ is said to be ${\rde C}$-acyclic if every morphism from $C$ to a complex in ${\rde C}$ is null-homotopic.}
\vskip 12pt
The two notions of locality are related by the following result:
\vskip 12pt
\noi{\bf 2.1 Lemma:} {\sl Let $C$ be an $A$-complex. Then $C$ is $W$-local if and only its homology is $W$-local.}
\vskip 12pt
\noi{\bf Proof:} For simplicity the module $\H_i(C)$ will be denoted by $H_i$. 

Let $K$ be a complex in ${\rde W}_0$ concentrated in degrees $p$ and $p+1$. Denote by $f: K_{p+1}\rightarrow K_p$ the non zero differential of $K$ and set:
$$C^K=\Hom(K,C)$$
We have a short exact sequence of complexes:
$$0\longrightarrow K_p\longrightarrow K\longrightarrow K_{p+1}\longrightarrow 0$$
and then a short exact sequence of graded differential modules:
$$0\longleftarrow\Hom(K_p,C)\longleftarrow C^K\longleftarrow\Hom(K_{p+1},C)\longleftarrow 0$$
inducing a long exact sequence:
$$\dots\rightarrow \H_i(C^K)\rightarrow\Hom(K_p,H_{p+i})\build\rightarrow_{}^{f^*}\Hom(K_{p+1},H_{p+i})\rightarrow\H_{i-1}(C^K)\rightarrow\dots$$
The graded module $C^K=\Hom(K,C)$ has trivial homology if and only if every morphism from $K$ to $C$ is null-homotopic. The desired result follows.\cqfd
\vskip 12pt
\noi{\bf Notations:} {\sl The class of $W$-local complexes in ${\rde C}_*(A)$ will be denoted by ${\rde L}$. The class of ${\rde L}$-acyclic complexes in ${\rde C}_*(A)$
will be denoted by ${\rde W}$. A morphism between two $A$-complexes will be called a ${\rde W}$-equivalence if its mapping cone belongs to ${\rde W}$.}
\vskip 12pt
\noi{\bf 2.2 Lemma:} {\sl The class ${\rde W}$ has the following properties:

$\bullet$ ${\rde W}$ contains ${\rde W}_0$ and the acyclic complexes.

$\bullet$ if $\ 0\longrightarrow C\longrightarrow C"\longrightarrow C''\longrightarrow 0$ is a short exact sequence of complexes, then if two of these complexes are in
${\rde W}$ the third one is in ${\rde W}$ too.

$\bullet$ If a complex $C$ is a direct sum of complexes in ${\rde W}$, then $C$ belongs to ${\rde W}$.

$\bullet$ If a complex $C$ is the colimit of an infinite sequence:
$$C_1\rightarrow C_2\rightarrow C_3\rightarrow \cdots $$
such that all these maps are cofibrations and all the complexes $C_i$ are in ${\rde W}$, then $C$ belongs to ${\rde W}$.}
\vskip 12pt
\noi{\bf Proof:} The first property is obvious. 

Consider a short exact sequence in ${\rde C}_*(A)$:
$$0\longrightarrow C\longrightarrow C'\longrightarrow C''\longrightarrow 0$$
Let $L$ be a $W$-local complex. Then we have a short exact sequence of graded differential modules:
$$0\longleftarrow\Hom(C,L)\longleftarrow\Hom(C',L)\longleftarrow\Hom(C'',L)\longleftarrow0$$
If two of the complexes are in ${\rde W}$, two of the graded differential modules $\Hom(C,L)$, $\Hom(C',L)$, $\Hom(C'',L)$ are acyclic. Then the third one is also
acyclic and the third complex lies in ${\rde W}$.

Suppose $C$ is a complex on the form:
$$C=\build\oplus_i^{} C_i$$
where each $C_i$ lies in ${\rde W}$. Since $C$ is a complex, there is en integer $p$ such that each $C_i$ vanishes in degree less than $p$. Let $L$ be a $W$-local
complex. We have:
$$\Hom(C,L)=\prod_i\Hom(C_i,L)$$
Since $C_i$ lies in ${\rde W}$, $\Hom(C_i,L)$ has trivial homology. Therefore the homology of $\Hom(C,L)$ is also trivial and every morphism from $C$ to $L$ is
null-homotopic. Hence $C$ belongs to ${\rde W}$.

Suppose $C$ is the colimit of an infinite sequence:
$$C_1\rightarrow C_2\rightarrow C_3\rightarrow \cdots $$
such that all these maps are cofibrations and each $C_i$ is in ${\rde W}$. Then there is an integer $p$ such that $C$ and each $C_i$ vanishes in degree
less than $p$. Let $L$ be a $W$-local complex. We have an isomorphism of graded differential modules:
$$\Hom(C,L)\simeq\ \build{\rm lim}_{\longleftarrow\atop k}^{}\Hom(C_k,L)$$
and then an exact sequence:
$$0\longrightarrow\ \build{\rm lim}_{\longleftarrow\atop k}^{}^1 \H_{i+1}(\Hom(C_k,L))\longrightarrow \H_i(\Hom(C,L))\longrightarrow\build{\rm lim}_{\longleftarrow\atop k}^{}
\H_i(\Hom(C_k,L))\longrightarrow 0$$
Since each $C_k$ lies in ${\rde W}$, each $\Hom(C_k,L)$ is acyclic and $\Hom(C,L)$ also. Therefore $C$ lies in ${\rde W}$.\cqfd
\vskip 12pt
\noi{\bf 2.3 Lemma:} {\sl The class ${\rde L}$ has the following properties:

$\bullet$ ${\rde L}$ contains the acyclic complexes.

$\bullet$ If $0\longrightarrow C\longrightarrow C'\longrightarrow C''\longrightarrow 0$ is a short exact sequence of complexes, then if two of these complexes are in
${\rde L}$, the third one is in ${\rde L}$ too.

$\bullet$ If a complex $C$ is a direct sum of complexes in ${\rde L}$, then $C$ belongs to ${\rde L}$.

$\bullet$  Suppose a complex $C$ is the colimit of an infinite sequence:
$$C_1\rightarrow C_2\rightarrow C_3\rightarrow \cdots $$
such that all these maps are cofibrations and each $C_i$ is in ${\rde L}$. Then $C$ belongs to ${\rde L}$.}
\vskip 12pt
\noi{\bf Proof:} The first property is obvious.

Consider a short exact sequence of complexes:
$$0\longrightarrow C\longrightarrow C'\longrightarrow C''\longrightarrow 0$$
If $K$ is a complex in ${\rde W}_0$, we have a short exact sequence:
$$0\longrightarrow\Hom(K,C)\longrightarrow\Hom(K,C')\longrightarrow\Hom(K,C'')\longrightarrow 0$$
If two of the complexes $C$, $C'$, $C''$ are in ${\rde L}$, two of the graded differential modules $\Hom(K,C)$, $\Hom(K,C')$, $\Hom(K,C'')$ are acyclic and the 
third one is also acyclic. Thus the third complex is in ${\rde L}$.

Suppose $C$ is the direct sum of complexes $C_i$ in ${\rde L}$. Let $K$ be a complex in ${\rde W}_0$. Since $K$ is finitely dominated, $\Hom(K,C)$ is the direct sum of
the graded differential modules $\Hom(K,C_i)$. Therefore $\Hom(K,C)$ is acyclic and $C$ belongs to ${\rde L}$.

Suppose $C$ is the colimit of the sequence: $C_1\rightarrow C_2\rightarrow\dots$. Let $K$ be a complex in ${\rde W}_0$. Since $K$ is finitely dominated, $\Hom(K,C)$ is
the colimit of the $\Hom(K,C_i)$ and $\H_*(\Hom(K,C))$ is the colimit of the $\H_*(\Hom(K,C_i))$. Therefore $\Hom(K,C)$ is acyclic and $C$ belongs to ${\rde L}$.\cqfd
\vskip 12pt
\noi{\bf Definition:} {\sl Let $C\rightarrow C'$ be a morphism in ${\rde C}_*(A )$. We say that $C'$ (or more precisely $C\rightarrow C'$) is a $W$-localization of
$C$ if the following holds:
 
$\bullet$ the morphism $C\rightarrow C'$ is a cofibration.

$\bullet$ $C'$ belongs to ${\rde L}$.

$\bullet$ $C\rightarrow C'$ is a ${\rde W}$-equivalence.}
\vskip12pt
\noi{\bf 2.4 Proposition:} {\sl Let $C$ be an $A$-complex and $\lambda_1:C\rightarrow C_1$ and $\lambda_2: C\rightarrow C_2$ be two $W$-localizations of $C$. Then there
exists a homotopy equivalence $f:C_1\rightarrow C_2$, unique up to homotopy, making the following diagram commutative:
$$\begin{tikzpicture}[scale=1/2] \draw (0,2) node {$C$} ; \draw (4,0) node {$C_2$} ; \draw (4,4) node {$C_1$} ; \draw [->,>=latex] (0.8,2.4)--(3.2,3.6) ;
\draw [->,>=latex] (0.8,1.6)--(3.2,0.4) ; \draw [->,>=latex] (4,3.2)--(4,0.8) ; \draw (2,3.5) node {$\lambda_1$}; \draw (2,0.5) node {$\lambda_2$};
\draw (4.5,2) node {$f$}; \end{tikzpicture}$$}
\vskip 12pt
\noi{\bf Proof:} Let $K_1$ be the Cokernel of $\lambda_1$. We have the following exact sequences:
$$0\longrightarrow C\longrightarrow C_1\longrightarrow K_1\longrightarrow0$$
$$0\longleftarrow\Hom(C,C_2)\longleftarrow\Hom(C_1,C_2)\longleftarrow\Hom(K_1,C_2)\longleftarrow0$$
Since $C_2$ is in ${\rde L}$ and $K_1$ in ${\rde W}$, the homology of $\Hom(K_1,C_2)$ is trivial and the restriction map:
$$\Hom(C_1,C_2)\longrightarrow \Hom(C,C_2)$$
is an isomorphism in homology. Therefore the morphism $\lambda_2$ can be lifted in a morphism $f:C_1\rightarrow C_2$. Moreover this lifting is unique up to homotopy.
Similarly, we produce a morphism $g: C_2\rightarrow C_1$. By unicity up to homotopy of $f\circ g$ and $g\circ f$, these two maps have to be homotopic to identities.
The result follows.\cqfd
\vskip 12pt
\noi{\bf 2.5 Theorem:} {\sl There exist a functor $\Phi$ from ${\rde C}_*(A)$ to itself and a morphism $\lambda$ from the identity of ${\rde C}_*(A)$ to $\Phi$
with the following properties:

$\bullet$ For every complex $C$, $\lambda_C:C\rightarrow \Phi(C)$ is a $W$-localization of $C$.

$\bullet$ The functor $\Phi$ sends cofibration to cofibration.

$\bullet$ If a complex $C$ vanishes in degree less than $p$, the same holds for $\Phi(C)$.}
\vskip 12pt
\noi{\bf Proof:} Let's take a set $X$ in ${\rde P}(A)$ such that, for every morphism $\alpha:P\rightarrow Q$ in $W$, $P$ and $Q$ are isomorphic
to some modules in $X$. Denote by $W'$ the set of morphisms $\alpha:P\rightarrow Q$ such that $P$ and $Q$ are in $X$ and there exist two isomorphisms $\varepsilon$ and
$\varepsilon'$ such that $\varepsilon\circ\alpha\circ\varepsilon'$ belongs to $W$.

The class of complexes of length $2$ such that the only non zero differential is in $W$ (resp. $W'$) will be denoted by ${\rde W}_0$ (resp. ${\rde W}'_0$). Notice that
the class ${\rde W}'_0$ is a set and every complex in ${\rde W}_0$ is isomorphic to a complex in ${\rde W}'_0$.

Let $C$ be a complex. Denote by $E(C)$ the set of pairs $(K,\varphi)$ where $K$ is a complex in ${\rde W}'_0$ and $\varphi$ is a non zero morphism of
degree $0$ from $K$ to $C$. The mapping cone of the map:
$$\oplus K\build\longrightarrow_{}^{\oplus\varphi}C$$
where the direct sum is over all $(K,\varphi)\in E(C)$, will be denoted by $\Phi_1(C)$ and the  complex $\oplus K$ will
be denoted by $\Psi(C)$. By construction we have an exact sequence of complexes:
$$0\longrightarrow C\build\longrightarrow_{}^i\Phi_1(C)\build\longrightarrow_{}^j\Psi(C)\longrightarrow0$$
where $i$ and $j$ are two morphisms of degree $0$ and $-1$ respectively.

Suppose $C$ vanishes in degree less than $p$. Then for every $(K,\varphi)\in E(C)$, $\varphi$ is not zero and the complex $K$ vanishes in degree less than $p-1$. Thus
$\Psi(C)$ vanishes in degree less that $p-1$ and $\Phi_1(C)$ vanishes in degree less than $p$. 

Let's denote by $[K,\varphi]\otimes K$ the component of $K$ in $\Psi(C)$ corresponding to $(K,\varphi)\in E(C)$.

If $f:C\rightarrow C'$ is a morphism, we get a morphism $f_*$ from $\Psi(C)$ to $\Psi(C')$ defined by:
$$\forall (K,\varphi)\in E(C),\ \forall x\in K,\ \ f_*([K,\varphi]\otimes x)=\left\{\matrix{[K,f\circ\varphi]\otimes x&\hbox{if}\ f\circ\varphi\not=0\cr 0&
\hbox{otherwise}\cr}\right.$$

Then $\Psi$ is a functor sending cofibration to cofibration. Moreover the map sending $[K,\varphi]\otimes x$ to $\varphi(x)$ is a morphism $\mu$ from $\Psi$ to the
identity of ${\rde C}_*(A)$ and $\Phi_1$ is the mapping cone of $\mu$. Thus $\Phi_1$ is a functor sending cofibration to cofibration.

The inclusion $C\subset\Phi_1(C)$ is a cofibration of degree $0$ and every morphism from a complex in ${\rde W}'_0$ to $C$ is null-homotopic in $\Phi_1(C)$. But every
complex in ${\rde W}_0$ is isomorphic to some complex in ${\rde W}'_0$. Hence every morphism from a complex in ${\rde W}_0$ to $C$ is null-homotopic in $\Phi_1(C)$.

By iterating this construction, we get an infinite sequence:
$$C\rightarrow\Phi_1(C)\rightarrow\Phi_2(C)\rightarrow\Phi_3(C)\rightarrow\cdots$$
where all these maps are cofibrations of degree $0$. Moreover each $\Phi_n(C)$ vanishes in degree less than $p$ and the colimit $\Phi(C)$ of this sequence is a well
defined complex in ${\rde C}_*(A)$ vanishing in degree less than $p$.

It is clear that $\Phi_1(C)/C$ is isomorphic to a direct sum of complexes in ${\rde W}_0$. Because of lemma 2.2, each $\Phi_{n+1}(C)/\Phi_n(C)$, each $\Phi_n(C)/C$ and
$\Phi(C)/C$ belongs to ${\rde W}$.

Let $K$ be a complex in ${\rde W}_0$ and $f$ be a morphism from $K$ to $\Phi(C)$. Since $K$ is finitely dominated $f(K)$ is contained in some $\Phi_n(C)$ and the map
$K\rightarrow \Phi_n(C)$ is null-homotopic in $\Phi_{n+1}(C)$. So every morphism from $K$ to $\Phi(C)$ is null-homotopic. Therefore $\Phi(C)$ is $W$-local.\cqfd
\vskip 12pt
\noi{\bf 2.6 Proposition:} {\sl The functor $\Phi$ is homotopically exact in the following sense:

$\bullet$ If a complex $C$ is the direct sum of complexes $C_i$, the morphism $\oplus_i \Phi(C_i)\rightarrow\Phi(C)$ is a homotopy equivalence.

$\bullet$ Let
$$0\longrightarrow C\build\longrightarrow_{}^i C'\build\longrightarrow_{}^j C''\longrightarrow0$$
be a short exact sequence of complexes. Then the morphism induced by $j$ from $\Phi(C')/i_*(\Phi(C))$ to $\Phi(C'')$ is a homotopy equivalence.}
\vskip 12pt
\noi{\bf Proof:} Suppose $C$ is the direct sum of complexes $C_i$. For each $i$ we have an exact sequence:
$$0\longrightarrow C_i\longrightarrow \Phi(C_i)\longrightarrow E_i\longrightarrow0$$
where $E_i$ lies in ${\rde W}$. So we get the following exact sequence:
$$0\longrightarrow C\longrightarrow \oplus_i\Phi(C_i)\longrightarrow \oplus_i E_i\longrightarrow0$$
By lemmas 2.2 and 2.3, $\oplus_i E_i$ belongs to ${\rde W}$ and $\oplus_i\Phi(C_i)$ belongs to ${\rde L}$. Therefore $\oplus_i\Phi(C_i)$ is a $W$-localization of $C$
and the map $\oplus\Phi(C_i)\rightarrow\Phi(C)$ is a homotopy equivalence.

Consider a short exact sequence:
$$0\longrightarrow C\build\longrightarrow_{}^i C'\build\longrightarrow_{}^j C''\longrightarrow0$$
The map $i$ induces a cofibration $i_*$ from $\Phi(C)$ to $\Phi(C')$. Denote by $\Sigma$ the quotient $\Phi(C')/i_*(\Phi(C))$. So we have the following commutative
diagram:
$$\diagram{0&\hfl{}{}&C&\hfl{i}{}&C'&\hfl{j}{}&C''&\hfl{}{}&0\cr
&&\vfl{}{}&&\vfl{}{\lambda'}&&\vfl{}{\lambda''}&&\cr
0&\hfl{}{}&\Phi(C)&\hfl{}{}&\Phi(C')&\hfl{}{}&\Sigma&\hfl{}{}&0\cr}$$
where the two lines are exact. Take a cofibration $\mu$ from $C''$ to a contractible complex $U$. So we get the following commutative diagram:
$$\diagram{0&\hfl{}{}&C&\hfl{i}{}&C'&\hfl{j}{}&C''&\hfl{}{}&0\cr
&&\vfl{}{}&&\vfl{}{\lambda'\oplus\mu\circ j}&&\vfl{}{\lambda''\oplus\mu}&&\cr
0&\hfl{}{}&\Phi(C)&\hfl{}{}&\Phi(C')\oplus U&\hfl{}{}&\Sigma\oplus U&\hfl{}{}&0\cr}$$
In this diagram, the two lines are exact and the vertical lines are cofibrations. Because of lemmas 2.2 and 2.3, $\Sigma\oplus U$ is local and the Cokernel of
$\lambda''\oplus\mu$ lies in ${\rde W}$. Therefore $C''\rightarrow \Sigma\oplus U$ is a $W$-localization of $C''$ and $\Sigma\rightarrow\Phi(C'')$ is
a homotopy equivalence.\cqfd
\vskip 12pt
\noi{\bf 2.7 Theorem:} {\sl There is a well defined graded ring $\Lambda_*$ and a ring homomorphism $A\rightarrow\Lambda_0$ such that:

$\bullet$ $\Lambda_i=0$ for $i<0$ and $A\rightarrow\Lambda_0$ is the Cohn localization of $A$ with respect to $W$.
  
$\bullet$ if $C$ is a $W$-localization of $A$ (considered as a complex concentrated in degree $0$), $\Lambda_*$ is isomorphic to $\H_*(C)$ as a graded right $A$-module.

$\bullet$ for every $W$-local complex $C$, the graded module $\H_*(C)$ is equipped with a structure of right $\Lambda_*$-module and every morphism 
$f:C\rightarrow C'$ between $W$-local complexes induces a $\Lambda_*$-linear map from $H_*(C)$ to $H_*(C')$.

$\bullet$ for every complex $C$, there is a canonical spectral sequence with the following $E^1$ and $E^2$ terms:
$$E_{pq}^1=C_p\build\otimes_A^{}\Lambda_q\hskip 36pt E_{pq}^2=\H_p(C\build\otimes_A^{} \Lambda_q)$$
converging to $H_{p+q}(\Phi(C))$. Moreover this spectral sequence is natural in $C$ and compatible with the $\Lambda_*$-action.

The graded ring $\Lambda_*$ is well defined up to canonical isomorphism and will be called the graded localization of $A$ with respect to $W$.}
\vskip 12pt
\noi{\bf Proof}: Let $C$ be a $W$-localization of $A$ (considered as a complex concentrated in degree $0$). Denote by $E$ the Cokernel of the cofibration
$A\rightarrow C$. We have two exact sequences:
$$0\longrightarrow A\longrightarrow C\longrightarrow E\longrightarrow0$$
$$0\longleftarrow\Hom(A,C)\longleftarrow\Hom(C,C)\longleftarrow\Hom(E,C)\longleftarrow0$$
But $C$ belongs to ${\rde L}$ and $E$ belongs to ${\rde W}$. Therefore $\Hom(E,C)$ has no homology and the composite map:
$$\Hom(C,C)\longrightarrow\Hom(A,C)\build\longrightarrow_{}^= C$$
is an isomorphism in homology. Because of Lemma 2.0, $\Hom(C,C)$ is a graded differential ring and $\H_*(\Hom(C,C))$ is a graded ring.

Suppose $C'$ is another localization of $A$. Then there exist two morphisms: $f:C\rightarrow C'$ and $g:C'\rightarrow C$ of degree $0$ compatible with the cofibrations
$A\rightarrow C$ and $A\rightarrow C'$. These morphisms are unique up to homotopy and $g$ is a homotopy inverse of $f$.
So we have a morphism $\varphi:\Hom(C,C)\rightarrow\Hom(C',C')$ defined by:
$$\varphi(u)=f\circ u\circ g$$
Since $f$ and $g$ are homotopy inverses to each other, $\varphi$ is a homotopy equivalence well defined up to homotopy. Moreover $\varphi$ is compatible, up to
homotopy, with the product and induces an isomorphism of rings between $\H_*(\Hom(C,C))$ and $\H_*(\Hom(C',C'))$. So the graded ring $\Lambda_*=\H_*(\Hom(C,C))$ is well
defined up to canonical isomorphism. 

Since $C$ vanishes in negative degrees, $\Lambda_i=\H_i(\Hom(C,C))\simeq\H_i(C)$ is trivial for $i<0$. 

Let $\lambda$ be the cofibration $A\rightarrow C$. Let $a$ and $b$ be two elements in $A$. These elements are two cycles in $A$ and $\lambda(a)$ and $\lambda(b)$ are
two cycles in $C$. These cycles define two morphisms from $A$ to $C$ sending $1$ to $\lambda(a)$ and $\lambda(b)$ and can be extended to two endomorphisms $\theta(a)$
and $\theta(b)$ of degree $0$ of $C$. So we have:
$$\theta(a)(\lambda(1))=\lambda(a)\hskip 24pt \theta(b)(\lambda(1))=\lambda(b)$$
and that implies:
$$\theta(a)\circ\theta(b)(\lambda(1))=\theta(a)(\lambda(b))=\theta(a)(\lambda(1)b)=\theta(a)(\lambda(1))b=\lambda(a)b=\lambda(ab)$$
Therefore $\theta(ab)$ is homotopic to $\theta(a)\circ\theta(b)$ and
the morphism from $A$ to $\H_0(C)\simeq\H_0(\Hom(C,C))=\Lambda_0$ is compatible with the ring structure. Thus $A\rightarrow\Lambda_0$ is a ring homomorphism.

Let $L$ be a local complex. We have the following exact sequence:
$$0\longleftarrow\Hom(A,L)\longleftarrow\Hom(C,L)\longleftarrow\Hom(E,L)\longleftarrow0$$
and, as above, the restriction map $\Hom(C,L)\rightarrow \Hom(A,L)=L$ is an isomorphism in homology. But $\Hom(C,C)$ acts on the right on $\Hom(C,L)$ and we get a
morphism from $\Hom(C,L)\otimes\Hom(C,C)$ to $\Hom(C,L)$ inducing on $\H_*(\Hom(C,L))=\H_*(L)$ a structure of right $\Lambda_*$-module.

As above one checks that this structure is well defined and every morphism between two local complexes induces in homology a $\Lambda_*$-linear morphism.  

Consider a complex $C$. Denote by $C_p$ the component of $C$ of degree $p$ and by $C(p)$ the $p$-skeleton of $C$. The complexes $\Phi(C(p))$ define a filtration
of $\Phi(C)$ and then a spectral sequence converging to the homology of $\Phi(C)$. The $E^1$ term of this spectral sequence is:
$$E^1_{pq}=\H_{p+q}(\Phi(C(p))/\Phi(C(p-1)))$$
Because or Proposition 2.6, we have:
$$E^1_{pq}\simeq\H_{p+q}(\Phi(C(p)/C(p-1)))=\H_{p+q}(\Phi(C_p))\simeq\H_{p+q}(C_p\build\otimes_A^{}\Phi(A))\simeq C_p\build\otimes_A^{}\Lambda_q$$
Moreover the differential $d^1$ is induced by the differential of $C$. So we have:
$$E^2_{pq}\simeq H_p(C\build\otimes_A^{}\Lambda_q)$$

Now the last thing to do is to prove that $\Lambda_0$ is the Cohn localization of $A$ with respect to $W$.

Let $R$ be a ring. If $M$ is a right $R$-module (resp. a left $R$-module), the dual module $D(M)=\Hom(M,R)$ is a left $R$-module (resp. a right $R$-module). The
correspondence $D$ is a contravariant functor and, if $M$ is finitely generated and projective, the module $D(M)$ is finitely generated and projective and the canonical
morphism $M\rightarrow D^2(M)$ is an isomorphism.

Let $f:A\rightarrow R$ be a ring homomorphism. If $M$ is a module in ${\rde P}(A)$, the right $R$-modules $D(\Hom(M,R))$ and $M\otimes_A R$ are
isomorphic. Therefore we have the following equivalences:
$$\forall (P\rightarrow Q)\in W,\ \ P\otimes_A R\build\longrightarrow_{}^\sim Q\otimes_A R$$
$$\Longleftrightarrow\ \ \forall (P\rightarrow Q)\in W,\ \ D(P\otimes_A R)\build\longleftarrow_{}^\sim D(Q\otimes_A R)$$
$$\Longleftrightarrow\ \ \forall (P\rightarrow Q)\in W,\ \ \Hom(P,R)\build\longleftarrow_{}^\sim \Hom(Q,R)$$
$$\Longleftrightarrow\ \ R\ \hbox{is}\ W\hbox{-local as a right}\ A\hbox{-module}$$
If this condition is satisfied consider a free resolution $C$ of the right $A$-module $R$. The complex $C$ is $W$-local and, because of the spectral sequence above, the
morphism $R\rightarrow R\otimes_A\Lambda_0$ is an isomorphism. Then we have:
$$\forall (P\rightarrow Q)\in W,\ \ P\otimes_A R\build\longrightarrow_{}^\sim Q\otimes_A R$$
$$\Longleftrightarrow\ \ \ R\build\longrightarrow_{}^\sim R\otimes_A\Lambda_0$$

In particular the morphism $\Lambda\rightarrow \Lambda\otimes_A\Lambda_0$ is an isomorphism and $\Lambda_0$ has the following property:
$$\forall (P\rightarrow Q)\in W,\ \ P\otimes_A \Lambda_0\build\longrightarrow_{}^\sim Q\otimes_A \Lambda_0$$
Hence the ring homomorphism $\lambda_0:A\rightarrow\Lambda_0$ factors through $\Lambda$:
$$A\build\longrightarrow_{}^\lambda\Lambda\build\longrightarrow_{}^\mu\Lambda_0$$
and we have a commutative diagram:
$$\begin{tikzpicture}[scale=1/2] \draw (0,4) node {$A$};\draw (4,4) node {$\Lambda$};\draw (0,0) node {$\Lambda_0$};\draw (4,0) node {$\Lambda_0$};
\draw [->,>=latex] (.8,4)--(3.2,4);\draw [->,>=latex] (.8,0)--(3.2,0);\draw [->,>=latex] (0,3.2)--(0,.8);\draw [->,>=latex] (4,3.2)--(4,.8);
\draw [->,>=latex] (.6,.6)--(3.4,3.4);\draw (2,4.4) node {$\lambda$};\draw (-.4,2) node {$\lambda_0$};\draw (4.4,2) node {$\mu$};\draw (2,.4) node {Id};
\draw (1.8,2.3) node {$f$};\end{tikzpicture}$$

In this diagram, all the objects are rings and $\Lambda$ and $\Lambda_0$ are right $\Lambda_0$-modules. Morphisms $\lambda$ and $\mu$ are ring homomorphisms and
$\mu$ is $\Lambda_0$-linear. The morphism $f$ is the composite:
$$\Lambda_0=A\otimes_A\Lambda_0\build\longrightarrow_{}^{\lambda\otimes1}\Lambda\otimes_A\Lambda_0\simeq\Lambda$$
and $f$ is $\Lambda_0$-linear.

Thus $\mu:\Lambda\rightarrow\Lambda_0$ is a split surjection and we have a decomposition of right $\Lambda_0$-modules: $\Lambda=\Lambda_0\oplus J$. With this
decomposition the maps $f$ and $\mu$ are given by:
$$\forall (a,u)\in\Lambda_0\times J,\ \ f(a)=a\oplus0,\ \ \mu(a\oplus u)=a$$

Since $\mu$ is a ring homomorphism the product on $\Lambda$ denoted by $\times$ can be written this way:
$$(a\oplus u)\times(b\oplus v)=(ab\oplus \delta(a,b)+a\times v+u\times b+u\times v)$$
with $a,b$ in $\Lambda_0$ and $u,v$ in $J$. On the other hand the multiplication on the left by an element $x\in \Lambda$ corresponds to a morphism between $W$-local
modules and this multiplication is a morphism of right $\Lambda_0$-modules. So we have:
$$(a\oplus0)\times(b\oplus 0)=(ab\oplus\delta(a,b))=(a\oplus0)\times(1\oplus 0)b=ab\oplus \delta(a,1)b\ \ \Longrightarrow\ \delta(a,b)=\delta(a,1)b$$
Since $\lambda$ is a ring homomorphism, $1\oplus0$ is equal to $1$ and that implies:
$$a\oplus 0=(a\oplus0)\times(1,0)=a\oplus\delta(a,1)\ \ \Longrightarrow\ \delta(a,1)=0$$
Hence the map $\delta$ is trivial and $f$ is a ring homomorphism.

Therefore $\Lambda_0\simeq f(\Lambda_0)$ is a subring of $\Lambda$ and every morphism in $W$ is sent to an isomorphism in ${\rde M}(\Lambda_0)$. By the universal
property, the ring $\Lambda_0$ is the Cohn localization $L(A,W)=\Lambda$ of $A$.\cqfd
\vskip 12pt
\noi{\bf 2.8 Corollary:} {\sl Let $A$ be a ring and $W$ be a class of morphisms in ${\rde M}(A)$. Let $\Lambda_*$ be the graded localization ring of $A$ with respect
to $W$ and $\Lambda=\Lambda_0$ be the corresponding Cohn localization. Let $p>0$ be an integer such that $\Lambda_i$
vanishes for $0<i<p$. Let $M$ be a $W$-local right $A$-module. Then we have isomorphisms of right $\Lambda$-modules:
$$M\simeq M\otimes_A \Lambda$$
$$\forall i, 0<i\leq p,\ \ \hbox{Tor}_i^A(M,\Lambda)\simeq 0$$
$$\hbox{Tor}_{p+1}^A(M,\Lambda)\simeq M\otimes_A\Lambda_p\simeq M\otimes_\Lambda\Lambda_p$$
Moreover, if $M$ is the module $\Lambda$, the isomorphism: Tor$_{p+1}^A(\Lambda,\Lambda)\simeq \Lambda_p$ is an isomorphism of $\Lambda$-bimodules.}
\vskip 12pt
\noi{\bf Proof:} Consider a projective resolution $C$ of the right $A$-module $M$. Since $C$ is local $\lambda_C:C\rightarrow\Phi(C)$ is a homotopy equivalence
and the spectral sequence of Theorem 2.7 implies the desired result.\cqfd
\vskip 12pt
\noi{\bf 2.9 Remark:} Because of this corollary, $\Lambda_i$ vanish for all $i>0$ if and only if $\Lambda$ is stably flat over $A$. An example of non stably flat
localization rings was found by Schofield [S] but theorem B implies a lot of explicit other examples. 
\vskip 12pt
A last result in this section is a description of classes ${\rde W}$ and ${\rde L}$ using only the ring $\Lambda=L(A,W)$:
\vskip 12pt  
\noi{\bf 2.10 Proposition:} {\sl The class ${\rde W}$ is the smallest class in ${\rde C}_*(A)$ satisfying the conditions of Lemma 2.2. Moreover for every $C$ in
${\rde C}_*(A)$, one has the following characterizations:
$$C\in{\rde W}\ \ \Longleftrightarrow\ \ \H_*(C\build\otimes_A^{}\Lambda)=0$$
$$C\in{\rde L}\ \ \Longleftrightarrow\ \ \H_*(C)\build\longrightarrow_{}^\sim H_*(C)\build\otimes_A^{}\Lambda$$}
\vskip 12pt
\noi{\bf Proof:} Let ${\rde W}'$ be the smallest class in ${\rde C}_*(A)$ satisfying the conditions of Lemma 2.2. Let $C$ be a complex in ${\rde W}$. By construction
the complex $\Phi(C)/C$ lies in ${\rde W}'$. Since $C$ belongs to ${\rde W}$, $\Phi(C)$ is contractible and belongs to ${\rde W}'$. Therefore $C$ belongs to ${\rde W}'$
too and ${\rde W}$ is the class ${\rde W}'$.

Let $C$ be a complex and $C'$ be the complex $C\otimes_A\Lambda$. Suppose $C'$ is not acyclic and let $p$ be the smallest integer such that $H_p(C')\not=0$. The spectral
sequence of Theorem 2.7 implies:
$$H_i(\Phi(C))\simeq\left\{\matrix{H_p(C')&\hbox{if}\ i=p\cr 0&\hbox{if}\ i<p\cr}\right.$$
If $C$ lies in ${\rde W}$, $\Phi(C)$ is acyclic and $H_*(C')$ has to be zero. Conversely, if $C'$ is acyclic, $\Phi(C)$ is acyclic too and $C$ belongs to ${\rde W}$.

Let $M$ be a local module and $C$ a free resolution of $M$. Because ot Lemma 2.1, $C$ lies in ${\rde L}$ and we have the following isomorphisms:
$$M\simeq \H_0(C)\simeq\H_0(\Phi(C))\simeq \H_0(C)\build\otimes_A^{}\Lambda\simeq M\build\otimes_A^{}\Lambda$$
If $C$ is a local complex, each homology module of $C$ is local and the map $\H_*(C)\rightarrow\H_*(C)\otimes_A\Lambda$ is an isomorphism.

Conversely, if $\H_*(C)$ is isomorphic to $\H_*(C)\otimes_A\Lambda$, each $\H_i(C)$ is local and $C$ belongs to ${\rde L}$.\cqfd

\vskip 24pt
\noi{\bf 3. Proof of the main theorems.}
\vskip 12pt
If $G$ is a group, the localization ring $L(\Z[G]\build\rightarrow_{}^\varepsilon\Z)$, where $\varepsilon$ is the augmentation map, will be denoted by $\Lambda(G)$.
\vskip 12pt
\noi{\bf 3.1 Lemma:} {\sl Let $G$ be a group and $H_1$ and $H_2$ be two finite subgroups of $G$ of coprime orders. Suppose $H_2$ normalizes
$H_1$. Then for every $x$ in $H_1$ and $y$ in $H_2$, $(1-x)(1-y)$ and $(1-y)(1-x)$ vanish in $\Lambda(G)$.}
\vskip 12pt
\noi{\bf Proof:} Let $n_1$ and $n_2$ be the orders of $H_1$ and $H_2$. Since these two numbers are coprime, there exist two integers $a$ and $b$ with: $an_1+bn_2=1$.
Consider the following elements in $\Z[G]$:
$$N_1=\build\sum_{z\in H_1}^{} z\hskip 24pt\omega_1=\build\sum_{z\in H_1}^{} (1-z)=n_1-N_1$$
$$N_2=\build\sum_{z\in H_2}^{} z\hskip 24pt\omega_2=\build\sum_{z\in H_2}^{} (1-z)=n_2-N_2$$
Since $y$ normalizes $H_1$, we have the following:
$$N_1x=xN_1=N_1\hskip 24pt N_2y=yN_2=N_2\hskip 24pt N_1y=yN_1$$
and then:
$$\omega_1(1-x)=(1-x)\omega_1=n_1(1-x)$$
$$\omega_2(1-y)=(1-y)\omega_2=n_2(1-y)$$
$$\omega_1y=y\omega_1$$
Let $U$ be the element $1-a\omega_1-b\omega_2\in\Z[G]$. This element is sent to $1$ in $\Z$ and becomes invertible in $\Lambda(G)$. Moreover $U$ commutes with $y$.
Thus we have:
$$(1-x)(1-y)U=(1-x)(1-a\omega_1-b\omega_2)(1-y)=(1-x)(1-an_1-bn_2)(1-y)=0$$
$$U(1-y)(1-x)=(1-y)(1-a\omega_1-b\omega_2)(1-x)=(1-y)(1-an_1-bn_2)(1-x)=0$$
Since $U$ is invertible in $\Lambda(G)$, $(1-x)(1-y)$ and $(1-y)(1-x)$ vanish in $\Lambda(G)$.\cqfd
\vskip 12pt
\noi{\bf 3.2 Lemma:} {\sl Let $G$ be a finite group. Then the image of $G$ in $\Lambda(G)$ is nilpotent. Moreover for every $x,y\in G$ of coprime orders, $(1-x)(1-y)$
vanishes in $\Lambda(G)$.}
\vskip 12pt
\noi{\bf Proof:} Let $G'$ be the image of the map $G\rightarrow\Lambda(G)$. Because of Proposition 1.11, the ring $\Lambda$ is the localization ring
$L(\Z[G']\rightarrow\Z)$. Because of Corollary 1.8, $G'$ doesn't contain any perfect subgroup and $G'$ is solvable. Thus there is a filtration:
$$1=H_0\subset H_1\subset\dots \subset H_n=G'$$
of $G'$ such that for each $k<n$, $H_k$ is a normal subgroup of $H_{k+1}$ and $H_{k+1}/H_k$ is a group of prime order.

Let $k$ be an integer with $0\leq k<n$. Suppose $H_k$ is nilpotent. Let $p$ be the order of $H_{k+1}/H_k$ and $S$ be a $p$-Sylow subgroup of $H_{k+1}$. The composite map
$S\rightarrow H_{k+1}\rightarrow H_{k+1}/H_k$ is surjective.

Let $q$ be any
prime. Since $H_k$ is nilpotent it has a unique $q$-Sylow subgroup $S_q$ and $H_k$ is isomorphic to the product of all the $S_q$'s. Actually the group $S_p$ is the group
$S\cap H_k$. Let $S'\subset H_k$ be the product of all the $S_q$'s with $q\not=p$, and $z$ be an element of $S$. Since $H_k$ is nilpotent, it has a unique $q$-Sylow
subgroup and the conjugation by $z$ sends $S_q$ to itself. Therefore $z$ lies in the normalizer of each $S_q$ and then in the normalizer of $S'$. 

Let $x$ be an element of $S'$. Because of Lemma 3.1, $(1-x)(1-z)$ and $(1-z)(1-x)$ vanish in $\Lambda(G')$ and $xz$ and $zx$ have the same image in $\Lambda(G')$. But
$G'$ is  included in $\Lambda(G')$ and we have: $xz=zx$. Therefore each element in $S$ commute with each element in $S'$. Let $H\subset H_{k+1}$ be the
product $S\times S'$. Since $H$ contains each $S_q$, it contains $H_k$. On the other hand $H$ contains the group $S$ which surjects onto $H_{k+1}/H_k$. Therefore $H$
is the group $H_{k+1}$ and $H_{k+1}=S\times S'$ is nilpotent. 

Thus each $H_k$ is nilpotent by induction and $G'$ is nilpotent too.

Let $x$ and $y$ be two elements in $G$ of coprime orders $a$ and $b$. Let $x'$ and $y'$ be the images of $x$ and $y$ in $G'$ and $H_1$ (resp. $H_2$) be the product of all
the $q$-Sylow subgroups of $G'$ with $q$ dividing $a$ (resp. $b$). Then $x'$ lies in $H_1$ and $y'$ lies in $H_2$. Since $G'$ is nilpotent, each element in $H_1$
commutes with each element in $H_2$ and, by Lemma 3.1, $(1-x')(1-y')$ vanishes in $\Lambda(G')=\Lambda(G)$. Thus $(1-x)(1-y)$ vanishes in $\Lambda(G)$.\cqfd

\vskip 12pt
Actually there is a more general setting of this result:
\vskip 12pt
\noi{\bf 3.3 Proposition:} {\sl Let $f:A\rightarrow B$ be a ring homomorphism and $\Lambda=L^h(A\rightarrow B)$ be the corresponding Cohn localization. Let $G$ be a
finite group contained in the group of units of $A$ and sent to $1$ in $B$. Then the image $G'$ of $G$ in $\Lambda$ is nilpotent and, for every $x,y\in G$ of coprime
orders, $(1-x)(1-y)$ vanishes in $\Lambda$. Moreover for each prime $p$ invertible in $B$, the order of $G'$ is coprime to $p$.}
\vskip 12pt
\noi{\bf Proof:} We have two commutative diagrams of rings:
$$\diagram{\Z[G]&\hfl{}{}&\Z\cr\vfl{}{}&&\vfl{}{}\cr A&\hfl{}{}&B\cr}\hskip 48pt\diagram{\Z[G]&\hfl{}{}&L(\Z[G]\rightarrow\Z)\cr\vfl{}{}&&\vfl{}{}\cr A&\hfl{}{}&
\Lambda\cr}$$
Since the image of $G$ in $L(\Z[G]\rightarrow\Z)$ is nilpotent by lemma 3.2, the image of $G$ in $\Lambda$ is also nilpotent. Moreover $(1-x)(1-y)$ vanishes
in $L(\Z[G]\rightarrow\Z)$ and then in $\Lambda$.

Suppose $p$ is a prime invertible in $B$. Let $S$ be a $p$-Sylow subgroup of $G'$. If $S$ is non trivial, there is an element $x\in G'$ of order $p$. We have:
$$0=1-x^p=(1-x)U$$
with: $U=1+x+x^2+\dots+x^{p-1}$. Then $U$ is send to $p$ in $B$ and $U$ is invertible in $\Lambda$. Thus $1-x$ vanishes in $\Lambda$ and $x$ is
equal to $1$. Hence $S$ is trivial and the order of $G'$ is coprime to $p$.\cqfd
\vskip 12pt
Let $G$ be a finite group. For any prime $p$, denote by $G_p$ the quotient of $G$ by the subgroup generated by the elements of $G$ of order coprime to $p$. The group
$G_p$ is the biggest $p$-group quotient of $G$. Denote by $G'$ the product of the $G_p$'s. The morphism $G\rightarrow G'$ is surjective and $G'$ is the universal
nilpotent quotient of $G$. Because of Proposition 1.11 and Lemma 3.2, the map $L(\Z[G]\rightarrow\Z)\rightarrow L(\Z[G']\rightarrow\Z)$ is an isomorphism. 

Denote by $A(G)$ the ring defined by the pull-back diagram:
$$\diagram{A(G)&\hfl{}{}&\prod_p\Z[G_p]\cr \vfl{}{}&&\vfl{}{\prod_p\varepsilon_p}\cr\Z&\hfl{\Delta}{}&\prod_p\Z}$$
where $\varepsilon_p$ is the augmentation map, $\Delta$ the diagonal map and the product is over all non trivial $G_p$. There is a canonical ring homomorphism from
$\Z[G]$ to $A(G)$. Moreover this morphism factors through $\Z[G']$. 
\vskip 12pt
\noi{\bf 3.4 Lemma:} {\sl The morphism $\Z[G]\rightarrow A(G)$ is onto and every element in its kernel is killed in $\Lambda(G)$.}
\vskip 12pt
\noi{\bf Proof:} Since $A(G]\rightarrow A(G')$ and $\Lambda(G)\rightarrow \Lambda(G')$ are isomorphisms, it is enough to prove the lemma when $G$ is nilpotent. In
this case $G$ is the product of the groups $G_p$.

For each group $\Gamma$ denote by $I(\Gamma)$ the augmentation ideal of $\Z[\Gamma]$. In order to prove that $\Z[G]\rightarrow A(G)$ is onto, it is enough to prove that
the induced morphism $I(G)\rightarrow\hbox{Ker}(A(G)\rightarrow\Z)$ is onto. But the kernel of $A(G)\rightarrow\Z$ is the product of the ideals $I(G_p)$.
So it is enough to show that the map $I(G)\rightarrow\prod_p I(G_p)$ is onto. 

Let $q$ be a prime and $u$ be an element in $I(G_q)$. The inclusion $G_q\subset G$ sends $u$ to an element $v\in I(G)$. Let $w=(w_p)$ be the image of $v$ in the 
product $\prod_p I(G_p)$. By construction, we have: $w_q=u$ and $w_p=0$ for every $p\not=q$. Therefore the map $I(G)\rightarrow\prod_p I(G_p)$ is onto.

For each prime $p$, denote by $I_p$ the image of $I(G_p)$ under the inclusion map $\Z[G_p]\subset \Z[G]$. Because of Lemma 3.1, $I_pI_q$ is killed in $\Lambda(G)$ for
every $p\not=q$. So denote by $J(G)$ the two-sided ideal of $\Z[G]$ generated by all the $I_pI_q$, $p\not=q$, and by $B(G)$ the quotient $\Z[G]/J(G)$. The ideal
$J(G)$ is killed in $\Lambda(G)$ and the map $\Z[G]\rightarrow\Lambda(G)$ factors through $B(G)$. By construction $J(G)$ is killed in $A(G)$ and the map $\Z[G]\rightarrow
A(G)$ factors through $B(G)$. So, in order to prove the lemma, it is enough to prove that the induced map $f:B(G)\rightarrow A(G)$ is an isomorphism.
\vskip 12pt
\noi{\bf Sub-lemma:} {\sl Let $x$ be an element of $G$. For each $p$ denote by $x_p$ the image of $x$ under the projection $G\rightarrow G_p$. Then $1-x-\sum_p(1-x_p)$
vanishes in $B(G)$.}
\vskip 12pt
Using this sub-lemma we are able to finish the proof of the lemma. Let $u\in\Z[G]$ be an element sent to $0$ in $A(G)$. Then $u$ is sent to $0$ in $\Z$ and $u$ lies in
$I(G)$. So we have a decomposition:
$$u=\sum_i\lambda_i(1-x_i)$$
where the $x_i$'s are in $G$ and the $\lambda_i$'s in $\Z$. Because of the sub-lemma, we may as well suppose that each $x_i$ lies in some $G_p$ and we have a
decomposition: $u=\sum_p u_p$, where each $u_p$ lies in $I(G_p)$.

Since $u$ vanishes in $A(G)$, it vanishes in $\Z[G_p]$ for each $p$ and $u_p$ is zero for each $p$. Therefore $u$ is sent to $0$ in $B(G)$ and $f$ is injective.
So $f$ is bijective and the lemma is proven.\cqfd
\vskip 12pt
\noi{\bf Proof of the sub-lemma:} Suppose $x$ is the product $x=y_1y_2\dots y_k$ where $y_i$ belongs to $G_{p_i}$, and all the primes $p_i$ are distinct. Denote by
$\equiv$ the equality in $B(G)$. We have:
$$1-x=1-y_1y_2\dots y_k=1-y_1y_2\dots y_{k-1}+y_1y_2\dots y_{k-1}(1-y_k)$$
But we have for every $i<k$: $y_i(1-y_k)\equiv 1-y_k$. Therefore we have:
$$1-x\equiv 1-y_1y_2\dots y_{k-1}+(1-y_k)$$
and, by induction, we get:
$$1-y_1y_2\dots y_k\equiv 1-y_1y_2\dots y_{k-1}+(1-y_k)\equiv1-y_1y_2\dots y_{k-2}+(1-y_{k-1})+(1-y_k)$$
$$\equiv \dots \equiv (1-y_1)+(1-y_2)+\dots+(1-y_k)$$
The result follows.\cqfd
\vskip 12pt
For each prime $p$, set:
$$n_p=\hbox{card}(G_p)\hskip 48pt N_p=\sum_{x\in G_p}x\hskip 48pt \omega_p=n_p-N_p=\sum_{x\in G_p}(1-x)$$
$\omega_p$ lies in the center of $\Z[G']$ and we have: $\omega_p^2=n_p\omega_p$. Consider the multiplicative sets $S_p=1+\Z\omega_p$ and 
$\Sigma_p=1+\Z n_p$. Since $n_p$ is a power of $p$, we have:
$$\Sigma_p^{-1}\Z=L^s(\Z\rightarrow\Z/n_p)=L^h(\Z\rightarrow\Z/n_p)=\Z_{(p)}$$
Let $R$ be the subring of $\Z[G']$ generated by all the $\omega_p$'s and $J$ be the ideal of $R$ generated by the $\omega_p$'s. The ring $R$ is contained in the center of
$\Z[G']$ and $R/J$ is the ring $\Z$. So we have another multiplicative set in the center of $\Z[G']$: the set $S=1+J$.

If $D$ is a commutative square in a category of modules, we say that $D$ is exact if $D$ is cartesian and cocartesian.
\vskip 12pt
\noi{\bf 3.5 Lemma:} {\sl Let $p$ be a prime. Then we have an exact square of right $\Z[G]$-modules:
$$\diagram{S^{-1}\Z[G_p]&\hfl{}{}&\Z_{(p)}[G_p]\cr\vfl{}{}&&\vfl{}{}\cr\Z&\hfl{}{}&\Z_{(p)}\cr}$$}
\vskip 12pt
\noi{\bf Proof:} The two-sided ideal of $\Z[G_p]$ generated by $N_p$ is the module $(N_p)=\Z N_p$.  Hence the following square is exact:
$$\diagram{\Z[G_p]&\hfl{}{}&\Z[G_p]/(N_p)\cr\vfl{}{}&&\vfl{}{}\cr\Z&\hfl{}{}&\Z/n_p\cr}$$
Since $\omega_p$ is killed in $\Z$, we have:
$$S_p^{-1}\Z=\Z\hskip 24pt S_p^{-1}\Z/n_p=\Z/n_p$$
But $S_p$ is equal to $\Sigma_p$ mod $N_p$. So we have:
$$S_p^{-1}(\Z[G_p]/(N_p))=\Sigma_p^{-1}(\Z[G_p]/(N_p))=\Z_{(p)}[G_p]/(N_p)$$
and we have an exact square:
$$\diagram{S_p^{-1}\Z[G_p]&\hfl{}{}&\Z_{(p)}[G_p]/(N_p)\cr\vfl{}{}&&\vfl{}{}\cr\Z&\hfl{}{}&\Z/n_p\cr}$$
Consider the following commutative diagram:
$$\diagram{S_p^{-1}\Z[G_p]&\hfl{}{}&\Z_{(p)}[G_p]&\hfl{}{}&\Z_{(p)}[G_p]/(N_p)\cr\vfl{}{}&&\vfl{}{}&&\vfl{}{}\cr\Z&\hfl{}{}&\Z_{(p)}&\hfl{}{}&\Z/n_p\cr}$$
In the square on the right the horizontal maps are surjective with isomorphic kernels. Then this square is exact. The composite square is also exact. Therefore the
square on the left is exact. On the other hand $\omega_q$ is killed in $\Z[G_p]$ for every $q\not=p$. Then $S_p^{-1}\Z[G_p]$ is equal to the ring $S^{-1}\Z[G_p]$ and the
lemma is proven.\cqfd
\vskip 12pt
\noi{\bf 3.6 Lemma:} {\sl There is an exact sequence of right $A(G)$-modules:
$$0\longrightarrow\build\oplus_p^{} I(G_p)_{(p)}\longrightarrow S^{-1}A(G)\longrightarrow \Z\longrightarrow0$$}
\vskip 12pt
\noi{\bf Proof:} We have the following exact sequence:
$$0\longrightarrow I(G_p)\longrightarrow\Z[G_p]\longrightarrow\Z\longrightarrow0$$
But Lemma 3.5 implies the exact sequence:
$$0\longrightarrow I(G_p)_{(p)}\longrightarrow S^{-1}\Z[G_p]\longrightarrow\Z\longrightarrow0$$
and $S^{-1}I(G_p)$ is isomorphic to $I(G_p)_{(p)}$.

By construction we have an exact sequence:
$$0\longrightarrow\build\oplus_p^{} I(G_p)\longrightarrow A(G)\longrightarrow \Z\longrightarrow0$$
and by inverting $S$ we get the desired exact sequence.\cqfd
\vskip 12pt
\noi{\bf Proof of Theorem A:}
\vskip 12pt
The pull-back diagram defining the ring $A(G)$ and the previous lemma imply that we have an exact square:
$$\diagram{S^{-1}A(G)&\hfl{}{}&\prod_p\Z_{(p)}[G_p]\cr\vfl{}{}&&\vfl{}{}\cr\Z&\hfl{}{}&\prod_p\Z_{(p)}\cr}$$
Thus the only thing to do is to prove that $S^{-1}A(G)$ is the Cohn localization of the augmentation map $\Z[G]\rightarrow\Z$.

Since $S$ may be seen as a multiplicative set in $A(G)$ sent to $1\in\Z$, we have the following:
$$L(\Z[G]\rightarrow\Z)\simeq L(\Z[G']\rightarrow\Z)\simeq L(A(G)\rightarrow\Z)\simeq L(S^{-1}A(G)\rightarrow\Z)$$

Recall that a ring homomorphism $f$ was called local if $f$ is surjective and every element in $f^{-1}(1)$ is invertible. Because of Proposition 1.10, each morphism
$\Z_{(p)}[G_p]\rightarrow\Z_{(p)}$ is local. Then the product of these maps is local and, because the square above is exact, the map $S^{-1}A(G)\rightarrow\Z$ is local
too. Then the desired result follows from Lemma 1.9.\cqfd
\vskip 12pt
It is also possible to determine the ring $L(\Z[G]\rightarrow\Z/p)$. Actually one has the following:
\vskip 12pt
\noi{\bf 3.7 Proposition:} {\sl Let $G$ be a finite group and $p$ be a prime. Then the localization ring $L(\Z[G]\rightarrow\Z_{(p)})=L(\Z[G]\rightarrow\Z/p)$ is
isomorphic to the ring $\Z_{(p)}[G_p]$.}
\vskip 12pt
\noi{\bf Proof:} Because of Proposition 3.3 and Lemma 1.9, we have:
$$L(\Z[G]\rightarrow\Z/p)\simeq L(\Z[G_p]\rightarrow\Z/p)$$
Since each integer coprime to $p$ is invertible in $\Z/p$, we have also:
$$L(\Z[G_p]\rightarrow\Z/p)\simeq L(\Z_{(p)}[G_p]\rightarrow\Z/p)\simeq L(\Z_{(p)}[G_p]\rightarrow\Z_{(p)}$$
and $L(\Z[G]\rightarrow\Z/p)$ is the ring $\Z_{(p)}[G_p]$ because of Proposition 1.10.\cqfd
\vskip 12pt
Actually it is possible to give a description of the graded localization ring of the augmentation morphism $\Z[G]\rightarrow\Z$. In order to do that we'll take the
following notations, for each group $G$, each prime $p$ and each pointed space $X$:

$\bullet$ $\Lambda(G)$ denote the localization ring of the augmentation morphism $\Z[G]\rightarrow\Z$ and $\Lambda(G)_*$ the corresponding graded localization ring.

$\bullet$ $\Lambda_{(p)}(G)$ denote the localization ring of the augmentation morphism $\Z[G]\rightarrow\Z_{(p)}$ and $\Lambda_{(p)}(G)_*$ the corresponding graded
localization ring.

$\bullet$ $BG$ denote the classifying space of $G$. This space is a pointed space.

$\bullet$ $\Omega X$ denote the loop space of $X$.

$\bullet$ $X_{(p)}$ denote the localization space of $X$ with respect to the homology $\H_*(-,\Z_{(p)})$.
\vskip 12pt
\noi{\bf 3.8 Theorem:} {\sl Let $G$ be a finite group. Then the graded localization ring $\Lambda(G)_*$ is determined by the following pull-back diagram of graded rings:
$$\diagram{\Lambda(G)_*&\hfl{}{}&\prod_p\Lambda_{(p)}(G)_*\cr\vfl{}{}&&\vfl{}{}\cr \Z&\hfl{}{}&\prod_p\Z_{(p)}\cr}$$
where the product is over all prime $p$ dividing the order of $G$.

Moreover, for each prime $p$, $\Lambda_{(p)}(G)_*$ is isomorphic to $\H_*(\Omega(BG_{(p)}),\Z_{(p)})$ as a graded ring.}
\vskip 12pt
\noi{\bf Proof:} Denote by ${\rde W}$ the class of $\Z[G]$-complexes $C$ such that $\H_*(C\otimes_{\Z[G]}\Lambda(G))$ vanishes and by ${\rde L}$ the corresponding class
of local complexes.

For each prime $p$, the class of $\Z[G]$-complexes $C$ such that $\H_*(C\otimes_{\Z[G]}\Lambda_{(p)}(G)$ vanishes will be denoted by ${\rde W}_p$ and the corresponding
class of local complexes will be denoted by ${\rde L}_p$. We have the following inclusions:
$${\rde W}\subset {\rde W}_p\hskip 48pt {\rde L}_p\subset {\rde L}$$

Denote by $H_p$ the kernel of the quotient map $G\rightarrow G_p$ and by $F_p$ the homotopy fiber of the map $BG\rightarrow BG_{(p)}$.
Since: $\H_1(H_p,\Z_{(p)})=0$, the localization map $BG\rightarrow BG_{(p)}$ is characterized by the following properties:
$$\pi_1(BG_{(p)})\simeq G_p$$
$$\forall i>1,\ \ \pi_i(BG_{(p)})\ \hbox{is a}\ \Z_{(p)}\hbox{--module}$$
$$\forall i>0,\ \ \H_i(F_p,\Z_{(p)})=0$$

Let $\widetilde{BG_{(p)}}$ be the universal cover of $BG_{(p)}$, we have the following:
$$\forall i>1,\ \ \pi_i(BG_{(p)})\ \hbox{is a}\ \Z_{(p)}\hbox{--module}$$
$$\Longleftrightarrow \forall i>1,\ \ \pi_i(\widetilde{BG_{(p)}})\ \hbox{is a}\ \Z_{(p)}\hbox{--module}$$
$$\Longleftrightarrow \forall i>0,\ \ \pi_i(\Omega(\widetilde{BG_{(p)}}))\ \hbox{is a}\ \Z_{(p)}\hbox{--module}$$
$$\Longleftrightarrow \forall i>0,\ \ \H_i(\Omega(\widetilde{BG_{(p)}}),\Z)\ \hbox{is a}\ \Z_{(p)}\hbox{--module}$$
$$\Longleftrightarrow \forall i>0,\ \ \H_i(\Omega(BG_{(p)}),\Z)\ \hbox{is a}\ \Z_{(p)}[G_p]\hbox{--module}$$
Let $\lambda$ be the map: $F_p\rightarrow BG$. Since $F_p$ is only defined up to homotopy, $\lambda$ may be chosen to be a fibration. Denote by $EG$ the universal cover
of $BG$. The group $G$ acts on the right on $EG$ with a free action and $BG$ is the quotient $EG/G$. Denote by $\widehat F_p$ the space defined by the pull-back diagram:
$$\diagram{\widehat F_p&\hfl{\widehat\lambda}{}&EG\cr\vfl{}{}&&\vfl{}{}\cr F_p&\hfl{\lambda}{}&BG\cr}$$

Actually, this diagram is a diagram of pointed soaces and the group $G$ acts freely on $\widehat F_p$. Then the pull-back in $\widehat F_p$ of the base point $x\in F_p$
is the discrete group $G$.

Let $M_p$ be a $\Z$-complex which is a free resolution of $\Z_{(p)}$. The inclusion $\Z\subset\Z_{(p)}$ and the product
on $\Z_{(p)}$ are induced by morphisms $\Z\rightarrow M_p$ and $M_p\otimes M_p\rightarrow M_p$.
Let $C_0$ be the chain complex $C_*(G,\Z)$ and $C_1$ be the chain complex $C_*(\widehat F_p,\Z)$. Because of the free $G$-action, $C_0$ and $C_1$ are
$\Z[G]$-complexes. Moreover the inclusion $\Z[G]\subset C_0$ is a homotopy equivalence. We will prove that the morphism: 
$$\Z[G]\subset C_0\subset C_1=C_1\otimes\Z\rightarrow C_1\otimes M_p$$
is a ${\rde W}_p$-localization.

We have:
$$H_*(C_1/C_0\build\otimes_{\Z[G]}^{} \Z_{(p)}[G_p])=\H_*(F_p,x;\Z_{(p)}[G_p])=\H_*(F_p,x;\Z_{(p)})\build\otimes_{\Z_{(p)}}^{}\Z_{(p)}[G_p]=0$$
and, because of Proposition 2.10, $C_0\rightarrow C_1$ is a ${\rde W}_p$-equivalence. But the morphism $C_1\simeq C_1\otimes\Z\rightarrow C_1\otimes M_p$ is a
${\rde W}_p$-equivalence. Then the morphism $\Z[G]\rightarrow C_1\otimes M_p$ is a ${\rde W}_p$-equivalence too.

On the other hand the fibers of $\lambda$ and $\widehat\lambda$ are both homotopically equivalent to the space $\Omega(BG_{(p)})$ and $\widehat F_p$ has the homotopy
type of $\Omega(BG_{(p)})$. Thus the homology $H_*(C_1\otimes M_p)\simeq\H_*(\Omega(BG_{(p)}),\Z_{(p)})$ is a $\Z_{(p)}[G_p]$-module.
Because of Proposition 2.10, $C_1\otimes M_p$ belongs to ${\rde L}_p$ and $C_1\otimes M_p$ is a ${\rde W}_p$-localization of $C_0$ and then of $\Z[G]$. Hence we have
isomorphisms of $\Z[G]$-modules:
$$H_*(C_1\otimes M_p)\simeq\Lambda_{(p)}(G)_*\simeq\H_*(\Omega(BG_{(p)}),\Z_{(p)})$$

Now we have to prove that the isomorphism $\Lambda_{(p)}(G)_*\simeq\H_*(\Omega(BG_{(p)}),\Z_{(p)})$ is an isomorphism of graded rings.

Let $C$ be the chain complex $C_*(\Omega(BG_{(p)}),\Z)$ and $C'$ be the complex $C\otimes M_p$. Then the morphism $\Z[G]\rightarrow C'$ is a ${\rde W}_p$-localization of
$\Z[G]$. On $C'$ there is two morphisms $\mu_1, \mu_2$ from $C'\otimes_\Z C'$ to $C'$: the morphism $\mu_1$ inducing on $\Lambda_{(p)}(G)_*$ the canonical structure of
graded ring and the morphism $\mu_2$ coming from the H-space structure of $\Omega(BG_{(p)})$ and the product $M_p\otimes M_p\rightarrow M_p$. Since the map
$G\rightarrow\Omega(BG_{(p)})$ is a morphism of H-spaces, the morphisms:
$$C'\build\otimes_\Z^{} \Z[G]\longrightarrow C'\build\otimes_\Z^{} C'\build\longrightarrow_{}^{\mu_i} C'$$
are both homotopic to the composite:
$$C'\build\otimes_\Z^{} \Z[G]\longrightarrow C'\build\otimes_{\Z[G]}^{}\Z[G]\simeq C'$$

Let $K$ be the cokernel of the morphism: $\Z[G]\longrightarrow C'$. This complex belongs to ${\rde W}_p$. The morphism $\mu_1-\mu_2$ is homotopic to a composite:
$$C'\build\otimes_\Z^{} C'\longrightarrow C'\build\otimes_\Z^{} K\build\longrightarrow_{}^\nu C'$$
But we have:
$$\Hom_{\Z[G]}(C'\build\otimes_\Z^{} K,C')\simeq\Hom_\Z(C',\Hom_{\Z[G]}(K,C'))$$
and $\nu$ can be seen as a morphism $\widetilde\nu:C'\longrightarrow\Hom_{\Z_{(p)}[G]}(K,C')$. But $K$ belongs to ${\rde W}_p$ and $C$ to ${\rde L}_p$. Then the graded
differential module $\Hom_{\Z[G]}(K,C')$ is acyclic. Hence $\widetilde\nu$ is a morphism from a graded differential free $\Z$-module to an acyclic graded
differential $\Z$-module and $\widetilde\nu$ and $\nu$ are homotopic to $0$. Therefore $\mu_1$ and $\mu_2$ are homotopic and $\Lambda_{(p)}(G)_*$ is isomorphic to
$\H_*(\Omega(BG_{(p)}),\Z_{(p)})$ as a graded ring.

If $p$ doesn't divide the order of $G$, $G_p$ is the trivial group and $\Omega(BG_{(p)})$ is contractible. In this case $\Lambda_{(p)}(G)_*$ is the graded ring
$\Z_{(p)}$.

Let ${\rde P}$ be the set of primes dividing the order of $G$. Denote by $C$ a ${\rde W}$-localization of $\Z[G]$ and, for each prime $p$ in ${\rde P}$, denote by $C_p$ a
${\rde W}_p$-localization of $\Z[G]$. The morphism $\Z[G]\rightarrow C_p$ is a cofibration.

If $E$ is a right $\Z[G]$-module, denote by $E_*$ a free resolution of $E$. Since $\H_0(C_p)$ is isomorphic to $\Z_{(p)}[G_p]$, the augmentation morphism
$\Z_{(p)}[G_p]\rightarrow\Z_{(p)}$ can be represented by a morphism $\varepsilon_p:C_p\rightarrow \Z_{(p)*}$. Since $C_p$ is only defined up to homotopy, we may as well
suppose that $\varepsilon_p$ is onto.

Let $\Delta$ be a morphism from $\Z_*$ to $\build\oplus_{p\in{\rde P}}^{}\Z_{(p)*}$ inducing on $\H_0$ the diagonal inclusion
$\Z\rightarrow\build\oplus _{p\in{\rde P}}^{}\Z_{(p)}$ and $C'$ be the complex defined by the pull-back diagram:
$$\diagram{C'&\hfl{}{}&\build\oplus_{p\in{\rde P}}^{}C_p\cr \vfl{}{}&&\vfl{}{\oplus\varepsilon_p}\cr \Z_*&\hfl{\Delta}{}&\build\oplus_{p\in{\rde P}}^{}\Z_{(p)*}\cr}$$

Since all the complexes $\Z_*$, $\Z_{(p)*}$, $C_p$ are in ${\rde L}$, $C'$ belongs to ${\rde L}$ too and we have a commurative diagram:
$$\diagram{\Z[G]&\hfl{}{}&C&\hfl{\lambda}{}&C'&\hfl{}{}&\build\oplus_{p\in{\rde P}}^{}C_p\cr&&&&\vfl{}{}&&\vfl{}{\oplus\varepsilon_p}\cr
&&&&\Z_*&\hfl{\Delta}{}&\build\oplus_{p\in{\rde P}}^{}\Z_{(p)*}\cr}$$

If $p$ is a prime, denote by $U_p$ a free $\Z$-resolution of $\Z/p$. For any $q\not=p$, we have:
$$\H_*(C_q\otimes U_p)\simeq \H_*(\Omega(BG_{(q)}),\Z_{(q)}\otimes\Z/p)=0$$
$$\H_*(\Z_{(q)}\otimes U_p)\simeq \Z_{(q)}\otimes\Z/p=0$$
and that implies that $C'\otimes U_p\rightarrow C_p\otimes U_p$ is a homotopy equivalence. Then we have:
$$\H_*(C'/\Z[G]\otimes U_p\build\otimes_{\Z[G]}^{}\Z)\simeq \H_*(C_p/\Z[G]\otimes U_p\build\otimes_{\Z[G]}^{}\Z)\simeq\H_*(C_p/\Z[G]\build\otimes_{\Z[G]}^{} U_p)$$
$$\simeq \H_*(C_p/\Z[G]\build\otimes_{\Z[G]}^{} (U_p)_{(p)})\simeq \H_*(C_p/\Z[G]\build\otimes_{\Z[G]}^{}\Z_{(p)}\otimes (U_p)_{(p)})=0$$
and $\lambda\otimes U_p$ is a ${\rde W}$-equivalence.

Let $K$ be the homotopy kernel of $\lambda$. The complex $K$ belongs to ${\rde L}$ and $K\otimes U_p$ belongs to ${\rde W}$ for every prime $p$. Therefore the
multiplication by a prime $p$ is a homotopy equivalence from $K$ to $K$ and $\H_*(K)$ is a $\Q[G]$-module.

Because of Theorem A, the morphism $\H_0(C)\rightarrow\H_0(C')$ is an isomorphism. Then the last thing to do is to prove that $H_i(C)$ and $H_i(C')$ are torsion
$\Z$-modules for all $i>0$.

In order to do that denote by $\Lambda$ the ring $\Lambda(G)$ or one of the rings $\Lambda_{(p)}(G)$. Consider the spectral sequence defined in Theorem 2.7 for a
resolution $C$ of $\Lambda\otimes\Q$. Since $\Q[G]$ is semisimple, this spectral sequence collapses and all the modules $\Lambda_i\otimes\Q$
vanish for $i>0$. Thus, for every $i>0$ and every prime $p$, $\H_i(C\otimes\Q)$ and $\H_i(C_p\otimes\Q)$ vanish. Hence  $H_i(C)$ and $H_i(C')$ are torsion
$\Z$-modules for all $i>0$ and the module $H_i(K)$ is torsion for every $i$ and $K$ is acyclic. Thus $C'$ is a ${\rde W}$-localization of
$\Z[G]$ and the theorem is proven.\cqfd

\vskip 12pt
\noi{\bf Proof of theorem B:}
\vskip 12pt
It is easy to check the implications:
$$\hbox{1)}\ \Longrightarrow\ \hbox{2)}\ \Longrightarrow\ \hbox{3)}\ \Longrightarrow\ \hbox{4)}$$

Suppose $G$ is nilpotent. Then the central elements $\omega_p$ are well defined in $\Z[G]$ and the multiplicative set $S$ also. Moreover $S$ lies in the center of
$\Z[G]$. Let $p$ and $q$ be two distinct primes and $x$ and $y$ be two elements in $G_p$ and $G_q$ respectively. We have the following:
$$\omega_p(1-x)=n_p(1-x)\hskip 38pt \omega_q(1-y)=n_q(1-y)$$
take two integers $a$ and $b$ with: $an_p+bn_q=1$. Then we have:
$$(1-a\omega_p-b\omega_q)(1-x)(1-y)=(1-an_p-bn_q)(1-x)(1-y)=0$$
and $(1-x)(1-y)$ is killed in $S^{-1}\Z[G]$. Therefore the ring $S^{-1}\Z[G]$ is the ring $S^{-1}(A(G))=\Lambda$. Thus the localization is a central localization
and 5) implies 1). Then the last thing to do is to prove: 4) implies 5).

Suppose $\Lambda$ is stably flat over $\Z[G]$. Consider the left $\Z[G]$-module $U=\Lambda\otimes_\Z\Lambda$ where the $G$-action on $U$ is defined by:
$$g(u\otimes v)=ug^{-1}\otimes g v$$
for every $(g,u,v)\in G\times\Lambda\times\Lambda$. An elementary computation shows that Tor$_i^{\Z[G]}(\Lambda,\Lambda)$ is isomorphic to $\H_i(G,U)$ and $\H_i(G,U)$
vanishes for all $i>0$.

For any prime $p$, one has also the $\Z[G]$-module: $V_p=\Z_{(p)}[G_p]\otimes_\Z\Z_{(p)}[G_p]$ where the $G$-action is defined by:
$$g(u\otimes v)=ug^{-1}\otimes g v$$
for every $(g,u,v)\in G\times\Z_{(p)}[G_p]\times\Z_{(p)}[G_p]$.

For each prime $p$, denote by $H_p$ the kernel of the quotient map $G\rightarrow G_p$. The intersection $H$ of the $H_p$'s is the kernel of the map 
$G\rightarrow G'=\prod_p G_p$.
Denote also by $E_p$ the ring $\Q[G_p]/(N_p)$, where $N_p$ is the sum of all elements in $G_p$.

Let $p$ be a prime. We have: $\H_i(G,U)_{(p)}=\H_i(G,U_{(p)})$ and $\H_i(G,U_{(p)})$ vanishes for $i>0$.  Because of theorem A, the ring homomorphism
$$\Lambda\longrightarrow\Z_{(p)}[G_p]\times\prod_{q\not=p}E_q$$
induces an isomorphism:
$$\Lambda_{(p)}\build\longrightarrow_{}^\sim\Z_{(p)}[G_p]\times\prod_{q\not=p}E_q$$
and the $\Z[G]$-module $U_{(p)}$ contains $V_p$ as a direct summand. Therefore $\H_i(G,V_p)$ vanishes for every $i>0$.

Let $W$ be the left $\Z[G]$-module $\Z_{(p)}[G_p]$ and $W_0$ be the module $W$ equipped with a trivial $G$-action. The morphism $V_p\longrightarrow W\otimes W_0$
defined by $u\otimes v\mapsto v\otimes uv$ for every $u$ and $v$ in $G_p$, is an isomorphism of $\Z[G]$-modules. Then we have for every $i>0$:
$$0\simeq\H_i(G,V_p)\simeq \H_i(G,W\otimes W_0)=\H_i(G,W)\otimes W_0$$
$$\Longrightarrow\ \ 0\simeq \H_i(G,W)\simeq\H_i(H_p,\Z_{(p)})$$
Because of the next lemma the order of $H_p$ is coprime to $p$. Therefore the order of $H$ is coprime to every prime and $H$ is the trivial group. Thus 
$G\simeq\prod_p G_p$ is nilpotent.\cqfd
\vskip 12pt
\noi{\bf Lemma:} {\sl Let $G$ be a finite group and $p$ be a prime. Suppose $\H_i(G,\Z_{(p)})=0$ for all $i>0$. Then the order of $G$ is coprime to $p$.}
\vskip 12pt
\noi{\bf Proof:} Since $\H_i(G,\Z_{(p)})$ vanishes for $i>0$, the Krull dimension of the $\F_p$-algebra $\H^*(G,\F_p)$ (or equivalently of its center) is zero.
But this dimension is known to be the maximal integer $n$ such that $G$ contains an elementary abelian $p$-group $(\Z/p)^n$ (see Quillen [Q] or Adem [A]).

Since this dimension is zero, $G$ doesn't contain any non-trivial elementary abelian $p$-group and the order of $G$ is coprime to $p$.\cqfd
\vskip 12pt
\noi{\bf Example 1:} Let $p$ be a prime and $n$ be an integer coprime to $p$. Let $G$ be a group of order $pn$ containing a cyclic subgroup $H$ of order $n$. Suppose
one has: $H=[G,H]$. Denote by $\Lambda$ the localization of the augmentation morphism $\Z[G]\rightarrow\Z$ and by $\Lambda_*$ the associated graded localization ring.
The group $G_p$ is the quotient $G/H$ and the ring $\Lambda$ is the fiber product:
$$\Lambda=\Z_{(p)}[G_p]\build\times_{\Z_{(p)}}^{}\Z$$
Moreover there is an element $\alpha$ of degree $2p-2$ such that:
$$\Lambda_*=\Lambda[\alpha]/(n\alpha)=\Lambda\oplus(\Z/n)\alpha\oplus(\Z/n)\alpha^2\oplus(\Z/n)\alpha^3\oplus\dots$$

In the spectral sequence of 2.7, the element $\alpha\in\Lambda_{2p-2}$ corresponds to a generator of the module:
$$\hbox{Tor}_{2p-1}^{\Z[G]}(\Lambda,\Lambda)\simeq\hbox{Ker}(\H_{2p-1}(G)\rightarrow\H_{2p-1}(G_p)\simeq\H_0(G_p,\H_{2p-1}(H))$$
$$\simeq\H_{2p-1}(H)\build\otimes_G^{}\Z\simeq\Z/n$$
\vskip12pt
\noi{\bf Example 2:} Let $G$ be the alternating group $A_5$ and $\widetilde G$ be its universal central extension. These groups are the unique perfect groups of order
$60$ and $120$. Denote by $\Lambda_*$ and $\widetilde\Lambda_*$ the graded localization rings of the augmentation maps $\Z[G]\rightarrow \Z$ and
$\Z[\widetilde G]\rightarrow\Z$. Then we have:
$$\widetilde\Lambda_*=\Z[x]/(120 x)$$
$$\Lambda_*=\Z[t,y]/(2t,t^2,60y)$$
with: $\partial^\circ t=1$, $\partial^\circ x=\partial^\circ y=2$. Moreover the morphism $\widetilde G\rightarrow G$ induces a morphism from $\widetilde\Lambda_*$ to
$\Lambda_*$ sending $x$ to $y$.
\vskip12pt
\noi{\bf Counter-example 1:} Let $H$ be a subgroup of a group $G$ such that $H=[G,H]$. If $G$ is finite, the group $H$ is killed in any nilpotent quotient of $G$ and
the localization of $\Z[G]\rightarrow\Z$ is isomorphic to the localization of $\Z[G/H]\rightarrow\Z$. But this property is not true in general even if $H$ is finite as
will be shown in the following example:

Let $p$ be an odd prime. Consider the group $G$ generated by two elements $x$ and $t$ with the only relations: $x^p=1$ and $tx=x^{-1}t$. The center $Z$ of $G$ is a free
group generated by $y=t^2$ and $G/Z$ is a dihedral group of order $2p$. Let $H\subset G$ be the subgroup generated by $x$. The group $H$ is isomorphic to $\Z/p$ and we
have: $H=[G,H]$. The Cohn localization of the augmentation map $\Z[G]\rightarrow\Z$ will be denoted by $\Lambda$. We will prove that the morphism $\Z[G]\rightarrow
\Lambda$ is injective and doesn't factors through $\Z[G/H]$.

Consider the field $F_1=\Q(t)$ of rational fractions in $t$. The correspondence $x\mapsto 1$ and $t\mapsto t$ induces a morphism $f_1$ from $\Z[G]$ to $F_1$. It is
clear that $F_1$ is a local $\Z[G]$-module.

Consider the field $F_2=\Q[\zeta_p](y)$ of rational fractions in $y$ with coefficients in the cyclotomic field $\Q[\zeta_p]$. The complex conjugation $\zeta_p\mapsto
\zeta_p^{-1}$ induces an involution $u\mapsto \overline u$ in $F_2$. Denote by $F_3$ the subfield of $F_2$ generated by $y$ and $\zeta_p+\zeta_p^{-1}$. The field $F_2$ is
a Galois extension of $F_3$ of degree $2$. Set: $A=F_2\oplus F_2 t$. This module is a left $F_2$-vector space. We define a product on $A$ by:
$$(a+bt)(c+dt)=ac+b\overline d y+(ad+b\overline c)t$$
and $A$ is a ring. The involution on $F_2$ extends to an anti-involution on $A$ by setting: $\overline t=-t$. It is easy to see that $u+\overline u$ and $u\overline u$
are in $F_3$ for each $u\in A$. Moreover for each non zero $u\in A$, $u\overline u$ is non zero too. Thus each non zero element $u\in A$ is invertible with inverse:
$$u^{-1}=\overline u(u\overline u)^{-1}$$
and $A$ is a division ring. Actually $A$ can be embedded in the skew field of quaternions. 

The correspondence $x\mapsto\zeta_p$ and $t\mapsto t$ induces a morphism $f_2$ from $\Z[G]$ to $A$. As above $A$ is a local $\Z[G]$-module.

Using $f_1$ and $f_2$ we have a morphism $f:\Z[G]\rightarrow F_1\times A$ which factors through $\Lambda$. But a straightforward computation shows that $f$ is injective
and the morphism $\Z[G]\rightarrow\Lambda$ is injective too.
\vskip12pt
\noi{\bf Counter-example 2:} Let $G$ be a group and $x$ and $y$ be two elements in $G$ with coprime orders. Denote by $\Lambda$ the localization ring of the
augmentation map $\Z[G]\rightarrow\Z$. If $G$ is finite, $x$ and $y$ commute in $\Lambda$ and $(1-x)(1-y)$ vanishes in
$\Lambda$. But this property is not true in general if $G$ is infinite.

Actually, if $G$ is defined by the following presentation: 
$$G=<x,y|x^2=1,y^3=1>$$
the elements $(1-x)(1-y)$ and $(1-y)(1-x)$ do not vanish in $\Lambda$ and $x$ and $y$ do not commute in $\Lambda$.

In order to prove that we'll consider the following rings:

$\bullet$ The ring $A=\Z+\Z j\subset\C$ generated by a third root of unit $j=\zeta_3$.

$\bullet$ The polynomial ring $B=A[t]$ and the ring $B'=S^{-1}B$ where $S$ is the multiplicative set $1+(1-t)B$. 

The ideal in $B$ generated by $(1+t,1-j)$ is denoted by $J$.

$\bullet$ The fraction field $K$ of $B$ (and of $B'$).

$\bullet$ The ring $R\subset M_2(K)$ defined by:
$$R=\Z\pmatrix{1&0\cr 0&1\cr}+\pmatrix{(1-t)B&B\cr (1-t^2)B&J\cr}$$
The map sending the unit matrix to $1$ and all the matrices in $\pmatrix{(1-t)B&B\cr (1-t^2)B&J\cr}$ to $0$, is a ring homomorphism $f:R\rightarrow\Z$.

The correspondence $x\mapsto \pmatrix{t&1\cr1-t^2&-t\cr}$ and $y\mapsto\pmatrix{1&0\cr0&j\cr}$ defines a ring homomorphism $\varphi:\Z[G]\rightarrow R$. Moreover the
morphism $f\circ\varphi$ is the augmentation map $\Z[G]\rightarrow\Z$. 

Denote by $\Lambda'$ the localization ring of $f:R\rightarrow\Z$. By naturality we have a commutative diagram:
$$\diagram{\Z[G]&\hfl{\varphi}{}&R&\hfl{}{}&M_2(K)\cr\vfl{}{}&&\vfl{}{}&&\cr\Lambda&\hfl{}{}&\Lambda'&&\cr}$$

Let $U$ be a $n\times n$ matrix with entries in $R$ sent by $f$ to the unit matrix in $M_n(\Z)$. This matrix may be seen as a matrix:
$$U=\pmatrix{\alpha&\beta\cr\gamma&\delta\cr}$$
with $\alpha$, $\beta$, $\gamma$, $\delta$ in $B_n=M_n(B)$ and:
$$\gamma\in(1-t^2)B_n\hskip 24pt \alpha-1\in(1-t)B_n\hskip 24pt \delta-1\in JB_n$$
Since $\alpha$ is congruent to $1$ mod $(1-t)B_n$, its determinant belongs to $S$ and $\alpha$ is invertible in $M_n(B')$. Then there exist a matrix
$\widehat\alpha\in M_n(B)$ and a polynomial $u\in B$ such that:
$$\alpha^{-1}=(1+(1-t)u)^{-1}\widehat\alpha$$
Thus we have in $B'$:
$$\det(U)=\det(\alpha)\det(\delta-\gamma\alpha^{-1}\beta)=\det(\alpha)(1+(1-t)u)^{-n}\det(\delta(1+(1-t)u)-\gamma\widehat\alpha\beta)$$
We have four ring homomorphisms: $g_+,g_-:B\rightarrow A$, $h_3:A\rightarrow \F_3$ and $h_4:A\rightarrow\F_4$ defined by:
$$g_+:t\mapsto 1\hskip 24pt g_-:t\mapsto -1\hskip 24pt h_3:j\mapsto 1\hskip 24pt h_4:2\mapsto0$$
Since $g_+$ sends $\alpha$ to $1\in$M$_n(A)$ and $1+(1-t)u$ to $1\in A$, $\det(\alpha)(1+(1-t)u)^{-n}$ is non zero in $B'$. We have:
$$g_-(\delta(1+(1-t)u)-\gamma\widehat\alpha\beta)=g_-(\delta(1+(1-t)u))=g_-(\delta)(1+2g_-(u))$$
$$\Longrightarrow\ \ g_-(\det(\delta(1+(1-t)u)-\gamma\widehat\alpha\beta))=g_-(\det(\delta))(1+2g_-(u))^n$$
On the other hand $1+2g_-(u)$ (resp. $g_-(\delta)$) is sent to $1$ by $h_4$ (resp. $h_3$). Therefore  $1+2g_-(u)$ and $g_-(\det(\delta))$ are non zero in $A$ and
$\det(\delta(1+(1-t)u)-\gamma\widehat\alpha\beta)$ is non zero in $B\subset B'$. 

Thus $\det(U)$ is non zero in $B'\subset K$ and $U$ is invertible in $M_{2n}(K)\simeq M_n(M_2(K))$. By the universal property, the ring homomorphism
$R\rightarrow M_2(K)$ factors through $\Lambda'$ and the ring homomorphism $\Z[G]\rightarrow M_2(K)$ factors through $\Lambda$.

A straightforward computation shows the following:
$$\varphi((1-x)(1-y))=(1-j)\pmatrix{0&-1\cr0&1+t\cr}$$
$$\varphi((1-y)(1-x))=(1-j)(1+t)\pmatrix{0&0\cr t-1&1\cr}$$
$$\varphi(xy-yx)=(1-j)\pmatrix{0&-1\cr 1-t^2&0\cr}$$
and, because these matrices are not zero in $M_2(K)$, the three elements $(1-x)(1-y)$, $(1-y)(1-x)$ and $xy-yx$ do not vanish in $\Lambda$.
\vskip 48pt
\noi{\bf References: }

\begin{list}{}{\leftmargin 55pt \labelsep 10pt \labelwidth 40pt \itemsep 0pt}
\item[{[A]}] A. Adem -- {\sl  Lectures on the Cohomology of Finite Groups} , Contemporary Mathematics {\bf 436} (2007), 317--334.
\item[{[Ba]}] H. Bass -- {\sl Algebraic K-theory}, W. A. Benjamin, coll. « Mathematics Lecture Note Series », (1968).
\item[{[Bo]}] A. K. Bousfield -- {\sl The localization of spaces with respect to homology}, Topology {\bf 14} (1975), 133--150.
\item[{[CR]}] S. Cappell and D. Ruberman -- {\sl Imbeddings and homology cobordisms of lens spaces}, Comm. Math. Helv. {\bf 63} (1988), 75--89.
\item[{[CS]}] S. Cappell and J. Shaneson -- {\sl The codimension two placement problem and homology equivalent manifolds}, Ann. of Math. {\bf 99} (1974), 277--348.
\item[{[C1]}] P. M. Cohn -- {\sl Free rings and their relations}, Academic Press (1985).
\item[{[C2]}] P. M. Cohn -- {\sl Inversive localization in Noetherian rings}, Comm. Pure. Appl. Math. {\bf 26} (1973), 679–691.
\item[{[D]}] W. G. Dwyer -- {\sl Noncommutative localization in homotopy theory}, Noncommutative localization in algebra and topology, LMS 
Lecture Note Series {\bf 330}, Cambridge University Press (2006), 24-39.
\item[{[FV]}] M. Farber and P. Vogel -- {\sl The Cohn localization of the free group ring}, Math. Proc. Camb. Phil. Soc {\bf 111} (1992), 433--443.
\item[{[HV]}] J. C. Hausmann and P. Vogel -- {\sl The plus construction and lifting maps from manifolds}, Proc. Symp. Pure Math. (AMS) {\bf 32} (1978), 67--76.
\item[{[LD]}] J. Y. Le Dimet -- {\sl Cobordisme d'enlacements de disques}, Bull. Soc. Math. France {\bf 116} (1988).
\item[{[NR]}] A. Neeman and A. Ranicki -- {\sl Noncommutative localization and chain complexes I. Algebraic K- and L-theory}, math.RA/0109118 (2001).
\item[{[Q]}] D. Quillen -- {\sl The spectrum of an equivariant cohomology ring I \& II}, Ann. Math. {\bf 94} (1971), 549--602.
\item[{[R]}] A. Ranicki -- {\sl Noncommutative localization in topology}, math.AT/0303046 (2003).
\item[{[S]}] A. H. Schofield -- {\sl Representation of rings over skew fields}, LMS Lecture Note Series {\bf 82}, Cambridge University Press (1985).
\item[{[T]}] P. Teichner -- {\sl Flatness and the Ore condition for rings}, Proceedings of the AMS, {\bf 131} (2003), 1977--1980.
\item[{[V1]}] P. Vogel -- {\sl Un th\' eor\` eme de Hurewicz homologique}, Comm. Math. Helv. {\bf 52}, (1977), 393--413.
\item[{[V2]}] P. Vogel -- {\sl On the obstruction group in homology surgery}, Publ. Math. I.H.E.S Vol 55 N$^\circ$ 1 (1982), 165--206.

\end{list}
\end{document}